# A Multidimensional Self-Adaptive Numerical Simulation Framework for Semiconductor Boltzmann Transport Equation

Zeyu Zhang; Xiaoyu Zhang*; Zhigang Song; Qing Fang

## Abstract

This research addresses the numerical simulation of the Boltzmann transport equation for semiconductor devices by proposing a multidimensional self-adaptive numerical simulation framework. This framework is applied to two important generalized forms of the equation: a parabolic equation with singular properties on the unit disk and a continuity equation. The study enhances the alignment of numerical simulations with physical characteristics through polar coordinate transformation and variable drift-diffusion coefficients. Innovatively, a multidimensional adaptive mesh partitioning strategy for radius-angle-time is designed and combined with an adjustable finite difference scheme to construct a highly adaptive numerical simulation method. In the construction of discrete schemes, the Swartztrauber-Sweet method and the control volume method are employed to effectively eliminate the origin singularity caused by polar coordinate transformation. On the programming front, a parallelized MATLAB algorithm is developed to optimize code execution efficiency. Numerical comparative experiments demonstrate that the adaptive method improves the accuracy of the parabolic equation by 1 to 7 times and that of the continuity equation by 10% to 70% while maintaining computational efficiency, significantly enhancing numerical simulation accuracy with high stability. Furthermore, this study systematically verifies the algorithm's convergence, stability, and parameter sensitivity using error visualization and other means. It also explores optimal parameters and establishes tuning optimization criteria. The research provides theoretical support for high-precision and highly adaptive methods in semiconductor device simulation, demonstrating outstanding advantages in handling singular regions.

# Content

# 1 Introduction

## 1.1 Research on the Background and Significance of Equation Objects

With the rapid advancement of microelectronics and semiconductor material technologies, microelectronic devices are playing an increasingly critical role in various sectors of the national economy, national defense construction, and daily life. However, the development process of traditional integrated circuits is complex and costly, requiring repeated experimentation to achieve satisfactory design outcomes. This has led to prolonged development cycles for semiconductor devices and significant challenges in optimization[1]. Under the development of "computer simulation of semiconductor devices," the introduction of computer-aided design (CAD) tools has significantly optimized this process, reducing both R&D time and costs while improving design accuracy. Since Shockley laid the theoretical foundation for semiconductor devices, analytical models have been widely used to analyze one-dimensional semiconductor devices. However, with the advancement of large-scale integration (LSI) and very-large-scale integration (VLSI) technologies, the reduction in device size has made two- and three-dimensional effects increasingly prominent. As a result, traditional analytical models are no longer adequate, and numerical simulation methods have gradually become a key technology for studying the characteristics of small-scale devices. The use of numerical simulation technology in semiconductor development not only shortens cycles and reduces costs but also helps in understanding the internal details of devices and discovering new physical phenomena[2].

This paper first examines the degenerate form of the Boltzmann transport equation: the parabolic differential equation, which holds significant importance in the numerical simulation of the semiconductor Boltzmann transport equation. It is also widely applied in fields such as physics, biology, atmospheric science, and environmental engineering. For example, this class of equations can be used to model temperature distribution in the thermosphere, gradient-driven models of tumor anti-angiogenesis, and air pollution modeling systems[3-5]. Therefore, the numerical solution of parabolic equations has always been a research hotspot.

Based on the research foundation of adaptive numerical simulation for parabolic differential equations, this paper investigates the continuity equation (drift-diffusion model) derived from the semiconductor Boltzmann transport equation on a unit disk with singular properties (at the boundary and within the disk). This derivation incorporates the relaxation time approximation, quasi-neutrality approximation, and low-field strength assumption. This model also holds significant importance in the numerical simulation of the Boltzmann transport equation. This field originated from the work of Shockley and Gummel[6,7]. Researchers developed the fundamental theory of PN junctions and, for the first time, fully described the mathematical model of semiconductor devices using partial differential equations[8]. The continuity equation (drift-diffusion model) is widely used in mathematical modeling and numerical simulation of semiconductor devices. It serves as one of the core frameworks for describing charge carrier (electrons and holes) transport behavior in semiconductors or the electrodynamic airflow processes in corona discharges. This model can be applied to simulate the internal mechanisms of devices such as diodes, transistors, and solar cells. In recent years, with the enhanced computational capabilities of research equipment, numerical simulation has gradually become a key tool for studying these complex systems[9]. Numerical simulation methods enable researchers to verify, analyze, and evaluate various functions of semiconductor devices within a relatively short time frame, making them critically important in semiconductor device modeling. Key research focuses and challenges in this field include handling degenerate parabolic equations in multidimensional spaces, analyzing the convergence and stability of numerical schemes, and designing efficient algorithms.

The study by Marciulionis et al. employed the finite difference method to approximate the Poisson equation in polar coordinates, establishing a model suitable for conventional personal computers. This model can simulate corona discharge characteristic curves and predict the velocity distribution of electrodynamic airflow. The computationally derived velocity values were further validated through experimental measurements[9].

This study also references other literature on numerical simulations of transport models based on continuity equations (drift-diffusion model) derived from the semiconductor Boltzmann transport equation[10-22].

## 1.2 Research Background on the Application of Finite Difference Methods to Equation Objects

The finite difference method holds a significant position in the field of numerical simulation. Early improvements to this method primarily focused on simplifying the model, such as approximating the motion of charge carriers using the drift-diffusion equations[23]. Researchers have been dedicated to improving the finite difference method in recent years to enhance simulation accuracy and efficiency. Su et al. proposed an improved consistent, conservative, non-oscillatory, and high-order finite difference scheme for simulating variable-density turbulent flows at low Mach numbers. This approach ensures computational precision and improves stability through specific discretization strategies[24]. Furthermore, researchers have also developed a high-order finite difference method for incompressible flows[25].

To better address complex physical problems, Karl et al. combined the finite difference method with other numerical techniques and developed a stochastic simulation-based approach for solving nonlinear drift-diffusion-recombination transport equations in semiconductors[26].

The accurate specification and treatment of boundary conditions are crucial when dealing with the boundaries of unit disks with singular properties. In their study, Takao et al. conducted detailed analysis and tailored treatments for different boundary scenarios, ensuring the accuracy of numerical simulation results[27]. Additionally, this paper cites other literature on the application of boundary conditions in the numerical solution of differential equations[28,29].

In the field of semiconductor finite difference simulations, Sridharan proposed a numerical simulation model based on a multi-layer finite difference method. This model can be used to analyze the power and ground plane characteristics of semiconductors in the presence of non-uniform dielectric and conductive regions. The approach enables relatively accurate simulations of electromagnetic behavior in complex semiconductor packaging structures[30]. Huang et al. utilized both the Plane Wave Expansion Method (PWEM) and the Finite Difference Method (FDM) to calculate the band structure of strained quantum well Semiconductor Optical Amplifiers (SOAs). They systematically compared the differences between these two methods in terms of accuracy, computational speed, and stability[31]. When studying Bloch oscillations in semiconductor superlattices, Degond et al. derived hydrodynamic equations describing the evolution of electron density by considering multiple factors, such as the electron distribution function, through the Boltzmann-Poisson transport model[32].

The scholars introduced an appropriate coordinate system and a non-dimensionalization method to make the equations more suitable for numerical simulation on a unit disk region. Additionally, they employed different finite difference schemes for various terms in the drift-diffusion equations to enhance numerical accuracy[33]. Mirzadeh et al. introduced certain improvements to the traditional finite difference scheme to accommodate complex boundary conditions and enhance the stability and accuracy of numerical computations[34] This paper also references other literature that applies finite difference numerical simulation to practical research[35,36].

## 1.3 Research Background on Adaptive Finite Difference Methods and Adaptive Meshes

In the finite difference method, the Crank-Nicolson scheme is a highly classical difference formulation. Researchers such as Hu et al. employed the finite difference method to study high-dimensional Caputo-type parabolic equations with fractional Laplace operators[37]. Tesfaye Aga Bullo et al. developed a parameter-uniform implicit scheme for singularly perturbed reaction-diffusion initial-boundary value problems in the spatial direction[38]. Mohammad Tamsir et al. proposed a Crank-Nicolson collocation method based on a modified quintic B-spline[39]. Bhal Santosh Kumar et al. developed a one-dimensional Crank-Nicolson orthogonal spline collocation method for interface problems[40]. Bondare A. S. et al. employed a method combining spatial discretization via arbitrary finite-dimensional subspaces in a separable Hilbert space with Crank-Nicolson temporal discretization to solve linear variational parabolic equations under periodic conditions[41]. Apel Thomas et al. investigated the Crank-Nicolson scheme for optimal control problems governed by evolutionary equations[42].

In recent years, many mathematicians, physicists, and researchers in other scientific fields have devoted efforts to exploring numerical methods for solving ordinary and partial differential equations with singular exact solutions. This study focuses on the initial-boundary value problem of the continuity equation (drift-diffusion model) derived from the semiconductor Boltzmann transport equation on the unit disk, in regions with singular properties. It is assumed that the derivatives tend to infinity at the boundary of the disk, and the interior region of the disk also exhibits certain singular characteristics. Zhang et al. investigated the case where the exact solution exhibits singular behavior at the disk boundary. Their results demonstrate that with appropriate local mesh refinement, the Swartztrauber-Sweet scheme can achieve nearly second-order accuracy[43]. Zha et al. studied the initial-boundary value problem for a one-dimensional system of semilinear wave equations satisfying the null condition. They proved that for homogeneous Dirichlet or Neumann boundary conditions and sufficiently small initial data, classical solutions always exist globally[44].

Researchers have made significant progress in finite difference methods (FDM) and superconvergence techniques for problems involving singular solutions with non-uniform spatial discretization and local mesh refinement. Literature studies by Yamamoto et al. indicate that coordinate transformations can enhance the accuracy of approximate solutions[45]. Additionally, Zhang et al. employed a stretching polynomial function with parameters to construct local mesh refinement. Their numerical experiments revealed the existence of an optimal parameter value that enables the approximate solution to achieve the best accuracy[46]. To overcome the initial singularity, Zhou et al. adopted the Alikhanov scheme based on a non-uniform temporal grid for time discretization[47]. For further studies on non-uniform spatial discretization methods and problems with singular solutions, please refer to other relevant literature[48-62].

## 1.4 Research Background on Polar Coordinate Transformation Equations

Polar coordinate transformation equations can help us better understand the variation behavior of physical quantities inside materials. When dealing with linear elasticity problems involving complex geometries or non-uniform media, this approach can describe the distribution of physical quantities within objects featuring circular or annular structures more accurately than the Cartesian coordinate system.

Vinh P. C. and Tung D. X. proposed a homogenization equation for linear elasticity problems based on polar coordinates, which is suitable for cases involving highly oscillating interfaces. By expressing the linear elasticity equations in matrix form and applying standard homogenization methods, they derived explicit homogenized equations along with their associated continuity conditions. This work holds significance for problems involving highly oscillating circular interfaces and contributes to the study of how material and interface parameters influence the coefficients of the homogenized equations[63].

Halliday et al. investigated the lattice Boltzmann model by introducing spatially and velocity-dependent microscopic terms to adjust macroscopic dynamics. They explored how to apply forcing strategies within the LB framework to achieve specific target forms. Their approach is applicable to simulations of flow in circular cross-section pipes and can be extended to other complex fluid dynamics problems[64].

Regarding the research on the transformation of equations into polar coordinates, this paper also references other relevant literature.[65,66].

## 1.5 Problems in the Current Research Field and the Objectives of This Study

### 1.5.1 Treatment of Singular Properties

In the unit disk region with singular properties, accurately describing and processing areas near singular points poses significant challenges. Traditional finite difference methods often fail to adequately handle singularities near the unit disk's singular positions, leading to deviations between computed physical quantities—such as electric fields and carrier densities—and actual behavior. Moreover, these methods predominantly rely on one-dimensional adaptive meshes.

This study hypothesizes that the exact solution to this problem exhibits certain singularities both at the disk boundary and within the disk itself. Under these conditions, conventional finite difference schemes become inconsistent, rendering classical error analysis inapplicable and necessitating specialized numerical approaches.

Departing from traditional one-dimensional non-uniform discretization, this research innovatively introduces a spatially adaptive mesh in both radial and angular dimensions. This dual-dimensional adaptive grid allows for denser node distribution near the disk boundary and regions with singular properties, thereby better capturing singular variations in these areas and significantly enhancing numerical accuracy.

### 1.5.2 Numerical Stability and Accuracy

During numerical simulation, the nonlinear and coupled terms in the drift-diffusion equations pose significant challenges to numerical stability, particularly in regions with singular properties where numerical oscillations or divergence are prone to occur. The complex form of the equations can further exacerbate this instability. For instance, when solving semiconductor problems under high electric fields, the nonlinear interaction between electric fields and carriers may lead to unstable numerical solutions. Although existing studies have employed certain high-order difference schemes, achieving higher accuracy remains difficult when dealing with singular regions in unit disks. Moreover, high-order schemes often increase computational complexity, necessitating a trade-off between accuracy and computational cost in practical applications.

Building upon traditional one-dimensional adaptive meshes, this study innovatively introduces spatially adaptive meshes (in both radial and angular dimensions), temporally adaptive meshes, and temporally adaptive differential iteration schemes. This multi-dimensional approach allows for denser mesh allocation in regions with severe numerical oscillations, enabling better capture of rapid temporal variations and singular properties within the disk. Consequently, it significantly enhances numerical accuracy while maintaining approximately the same computational cost.

**1.5.3 Computational Efficiency**

With increasing demands for accuracy in semiconductor device simulations, finite difference methods require finer mesh partitioning—particularly near regions with singular properties—leading to a substantial increase in computational workload. When handling large-scale problems, the efficiency of matrix operations and other computational processes in finite difference numerical simulations is relatively low. Factors such as matrix size and sparsity can adversely affect computational speed, and traditional serial computing approaches struggle to meet the demands of large-scale numerical simulations. Thus, improving computational efficiency remains an urgent challenge.

The multi-dimensional adaptive method proposed in this study allows for adjustable mesh density across dimensions, precisely targeting regions where dense meshing is critical for accuracy. This approach significantly enhances numerical accuracy without increasing the total number of grid points. Furthermore, the MATLAB implementation in this study utilizes multi-threaded parallel computing code, substantially improving operational efficiency.

**1.5.4 Treatment of Cartesian-to-Polar Coordinate Transformation**

Traditional numerical simulations of the Boltzmann transport equation are primarily conducted in Cartesian coordinates. However, for symmetric physical systems such as circular or annular geometries (e.g., ring-shaped electrodes, cylindrical semiconductor devices), transforming the equation into polar coordinates can more naturally align the numerical simulation with the geometric structure, reduce discretization errors, simplify mesh generation, and improve both accuracy and interpretability. Under axisymmetric conditions, the angular derivative term vanishes, reducing the equation to a lower-dimensional problem involving only the radial variable, which significantly decreases computational complexity.

However, after the polar coordinate transformation, the term involving $r$ in the denominator introduces numerical singularity at the origin ($r=0$), which may lead to numerical instability. Despite this challenge, there is limited research on numerical solutions for continuity equations (e.g., drift-diffusion models) in polar coordinates. This gap necessitates further development in this area.

To address these issues, this study develops a novel multi-dimensional adaptive mesh method for such equations. This approach enables mesh refinement near the origin to capture rapidly varying physical quantities. Additionally, the Swartztrauber-Sweet method and an integral approach (finite volume method) are employed to construct difference schemes at the origin, effectively eliminating singularities and enhancing the accuracy of numerical simulations.

**1.5.5 Adaptive Variable Treatment of Drift and Diffusion Coefficients**

Traditional approaches often treat drift and diffusion coefficients as fixed values or set them piecewise based on empirical formulas, rather than considering them as variables that depend on spatial position. Treating these coefficients as functional variables significantly increases the complexity of the problem-solving process. However, it also greatly enhances the algorithm's adaptability to semiconductor devices and improves simulation accuracy. This study focuses on developing numerical algorithms that incorporate variable forms of drift and diffusion coefficients.

**1.5.6 Temporally Adaptive Differential Iteration Methods**

The temporally adaptive differential iteration method in this study is based on the Crank-Nicolson method. Existing literature offers limited research on optimizing the weight coefficients of this method, and integrating the Crank-Nicolson approach with multi-dimensional adaptive non-uniform grids to address initial-value problems with singularities introduces new complexities. This study builds upon the Crank-Nicolson method by incorporating an adaptive difference scheme with adjustable weight coefficients. It further investigates convergence analysis, exploration of optimal weight coefficients, and the synergistic effects of combining this method with adaptive non-uniform grids for singular initial-value problems, aiming to maximize the overall effectiveness of the approach.

## 1.6 Technical Route

This paper applies the multidimensional adaptive numerical simulation framework to two important extended forms of the Boltzmann transport equation to test the algorithm's effectiveness, with the technical approach illustrated in Figure 1.1.

## 1.7 Research Content

This paper aims to address the challenges in the numerical solution of the semiconductor Boltzmann equation, including the treatment of singularities, numerical stability and accuracy, computational efficiency, generalization capability of numerical simulations, and the handling of practical scenarios in simulations. To tackle these issues, this study innovatively proposes a multidimensional adaptive numerical simulation framework. The framework primarily includes: multidimensional adaptive grids in radius, angle, and time dimensions; a time-adaptive finite difference iterative method; parallel computing programs; polar coordinate transformation of the equation; variable adaptive drift-diffusion coefficients; and the use of Swarztrauber-Sweet and control volume methods to eliminate origin singularity. The algorithm is applied to two important generalized forms of the semiconductor Boltzmann transport equation. Through systematic convergence analysis and numerical experiments (including result visualization, comparative analysis, and parameter sensitivity analysis), the superiority of the algorithm is validated, and a basis for parameter tuning is provided. In summary, the multidimensional adaptive numerical simulation technique presented in this paper contributes to the development of highly adaptive, high-precision, high-stability, and high-efficiency numerical simulation technologies.

## 1.8 Thesis Structure Outline

In Chapter 1, this paper systematically elaborates on the background and significance of the studied equations (the Boltzmann equation, parabolic equation, and continuity equation). It summarizes the current state of domestic and international research on finite difference methods, adaptive finite difference methods, adaptive grids, and polar coordinate transformation methods in the context of these equations. Additionally, it identifies existing problems in current research, thereby introducing the research objectives and content of this paper.

In Chapter 2, this paper systematically explains the multidimensional adaptive solving process for semiconductor parabolic differential equations with singular properties on the boundary of a unit disk. The process primarily includes polar coordinate transformation of the equation, adaptive grid generation, derivation of discrete schemes, convergence analysis, and numerical experiments and analysis.

In Chapter 3, this paper aims to address the multidimensional adaptive numerical simulation of semiconductor continuity equations derived from the Boltzmann equation with singular properties on a unit disk. The process mainly involves polar coordinate transformation and simplification of the equation, adaptive grid generation, derivation of discrete schemes, establishment of iterative equation systems, and numerical experiments and analysis.

In Chapter 4, this paper provides a systematic summary of the research conducted in the first three chapters and offers insights into the potential applications of the multidimensional adaptive numerical simulation framework.

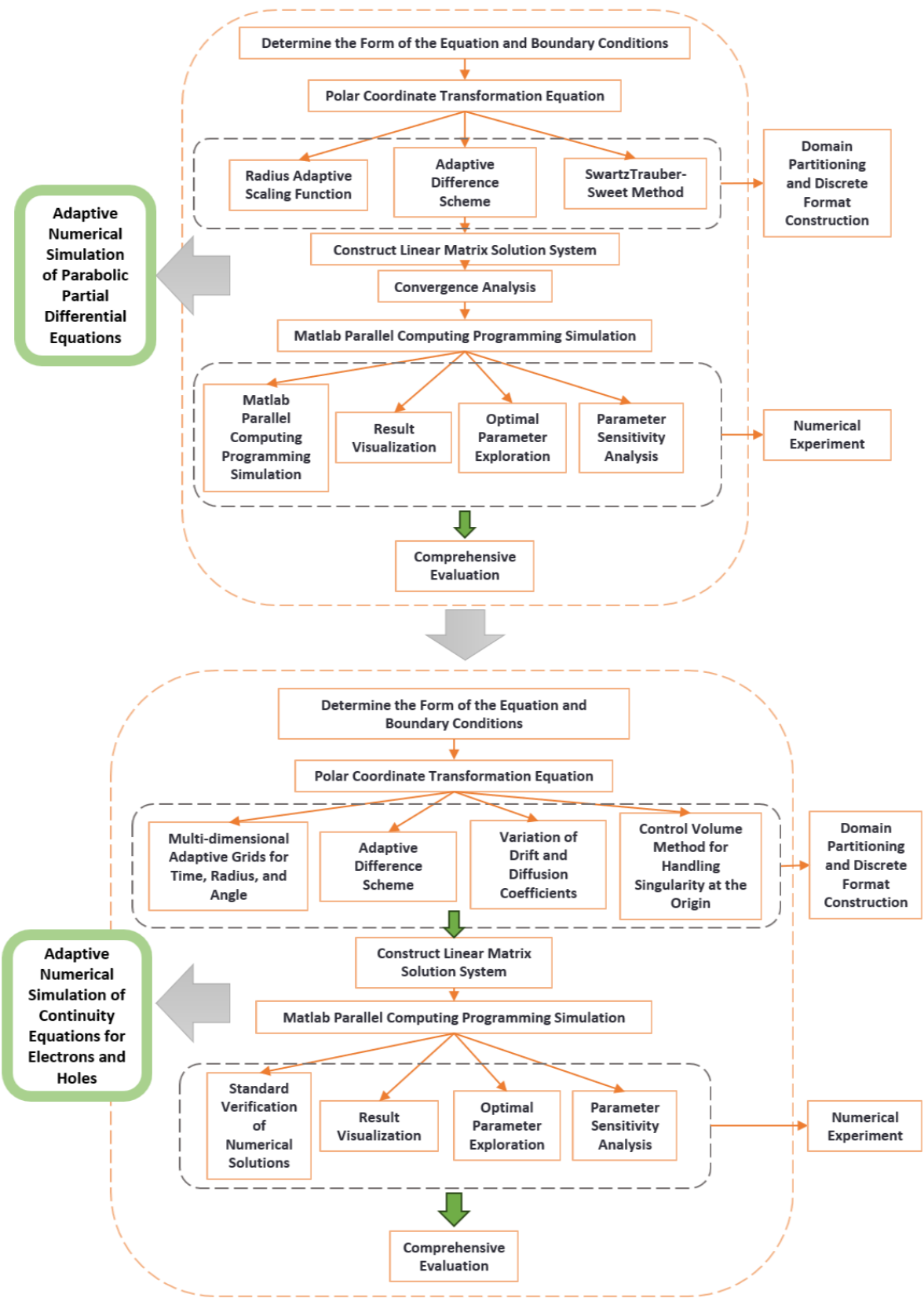

**Figure 1.1 Technical route**

# 2 Adaptive Numerical Simulation of Parabolic Differential Equations on the Unit Disk

This chapter systematically elucidates the processing of the multidimensional adaptive numerical simulation framework for parabolic equations in this paper. A series of numerical experiments are conducted to validate the superiority of the algorithm, along with certain analyses and discussions on its effectiveness and parameter properties.

## 2.1 Equation Form, Mesh Generation and Discretization Scheme

This chapter investigates the initial-boundary value problem for a parabolic equation on the following disk:

$$\frac{\partial u}{\partial t} - \Delta u + \tilde{b}(x,y)u = \tilde{f}(t,x,y), \quad 0 < t < T \in \tilde{\Omega},$$

$$u(t,x,y) = \tilde{\gamma}(t,x,y), \quad 0 < t < T, \quad (x,y) \in \partial\tilde{\Omega},$$

$$u(0,x,y) = u_0(x,y), \quad (x,y) \in \tilde{\Omega}.$$

The above equations satisfy $\tilde{b}(x,y) \geq 0$ in domain $\tilde{\Omega}$. $\tilde{f}$ and $\tilde{\gamma}$ are given functions. $\tilde{\Omega} = \{(x,y) \mid x^2 + y^2 < 1\}$.

Using the polar coordinate transformation method for equations, the above equation can be rewritten as:

$$\frac{\partial u}{\partial t} - \left[\frac{1}{r}\frac{\partial}{\partial r}\left(r\frac{\partial u}{\partial r}\right) + \frac{1}{r^2}\frac{\partial^2 u}{\partial \theta^2}\right] + b(r,\theta)u = f(t,r,\theta), \quad 0 < t < T, \quad (r,\theta) \in \Omega, \quad (2.1)$$

$$u(t,\tilde{r},\theta) = \gamma(t,\theta), \quad 0 < t < T, \quad (\tilde{r},\theta) \in \Gamma, \tag{2.2}$$

$$u(0,r,\theta) = u_0(r,\theta), \quad (r,\theta) \in \Omega, \tag{2.3}$$

where $\Omega = \{(r,\theta) \mid 0 < r < 1, 0 \leq \theta < 2\pi\}$ and $\Gamma = \{(\tilde{r},\theta) \mid \tilde{r} = 1, 0 \leq \theta < 2\pi\}$.

The radial expansion equation is written in the following form:

$$\upsilon(s) = \sin(\frac{s^p \pi}{2}), \quad 0 \leq s \leq 1, \quad p > 0, \tag{2.4}$$

under the following conditions: $\upsilon(0) = 0$, $\upsilon(1) = 1$。

The region $\Omega$ is then divided and processed as follows:

$$h = \frac{1}{m+1}, \quad r_i = \upsilon(ih), \quad r_{i+\frac{1}{2}} = \frac{r_i + r_{i+1}}{2}, \quad i = 0,1,2,\cdots,m+1,$$

$$h_i = r_i - r_{i-1}, \quad i = 1,2,\cdots,m+1,$$

$$\mu = \frac{2\pi}{n}, \quad \theta_j = j\mu, \quad j = 0,1,2,\cdots,n,$$

$$\Delta t = \frac{1}{k}, \quad k = 0,1,2,\cdots,K, \quad K \in Z^+,$$

and the divided spatial region is shown in Figure 2.1.

For the direction $t$, using an adaptive Crank-Nicolson method with an adjustable coefficient, and applying the Swartztrauber-Sweet method to equations (2.1)-(2.2), the following expressions are obtained:

$$\frac{U_{i,j}^{k+1} - U_{i,j}^k}{\Delta t} + a\wp_h U_{i,j}^k + (1-a)\wp_h U_{i,j}^{k+1} = af_{i,j}^k + (1-a)f_{i,j}^{k+1}, \quad k \in Z^+, \ 0 < a < 1, \quad (2.5)$$

where,

$$\wp_h U_{i,j} = \left( \frac{2}{r_i(h_i + h_{i+1})} \left( \frac{r_{i+\frac{1}{2}}}{h_{i+1}} + \frac{r_{i-\frac{1}{2}}}{h_i} \right) + \frac{2}{r_i^2 \mu^2} + b_{i,j} \right) U_{i,j} - \frac{2r_{i-\frac{1}{2}}}{r_i h_i(h_i + h_{i+1})} U_{i-1,j} - \frac{2r_{i+\frac{1}{2}}}{r_i h_{i+1}(h_i + h_{i+1})} U_{i+1,j}$$

$$- \frac{1}{r_i^2 \mu^2} U_{i,j-1} - \frac{1}{r_i^2 \mu^2} U_{i,j+1}, \quad i = 1,2,\cdots,m;\ j = 0,1,2,\cdots,n-1,$$

$$\wp_h U_{0,0} = \left( \frac{4}{h_1^2} + b_{0,0} \right) U_{0,0} - \frac{4}{nh_1^2} \sum_{j=0}^{n-1} U_{1,j},$$

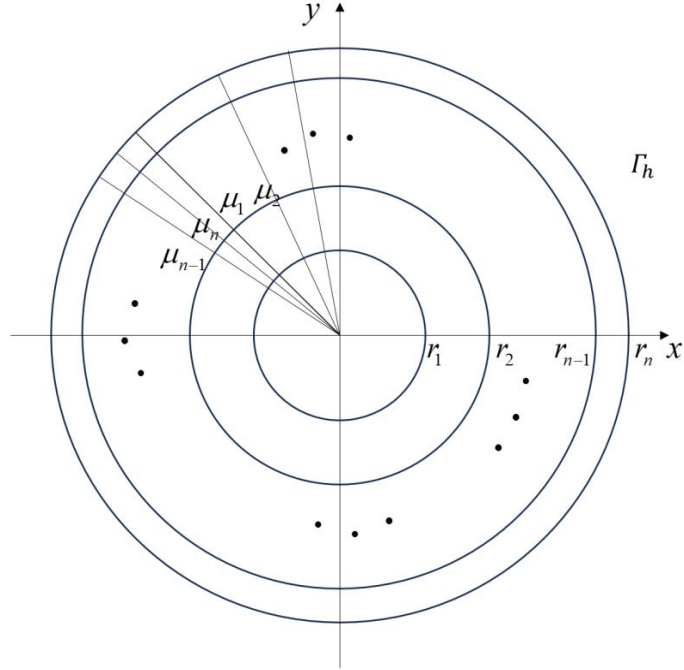

**Figure 2.1 Grid partitioning plot**

Moreover, $U_{i,j}^k$ is the approximate value of the exact solution $u(t_k, r_i, \theta_j)$ of equations (2.1)-(2.2) at time level $t_k = k\Delta t$, and $f_{i,j}^k = f(t_k, r_i, \theta_j)$. The above discrete equation satisfies the following initial and boundary conditions:

$$U_{i,n}^k = U_{i,0}^k, \tag{2.6}$$

$$U_{i,-1}^k = U_{i,n-1}^k, \quad i = 0,1,2,\cdots,m+1, \tag{2.7}$$

$$U_{0,j}^k = U_{0,0}^k, \tag{2.8}$$

$$U_{m+1,j}^k = \gamma_j^k, \quad j = 0,1,2,\cdots,n, \tag{2.9}$$

$$U_{i,j}^0 = u_0(r_i, \theta_j), \quad i = 0, 1, 2, \cdots, m+1; j = 0, 1, 2, \cdots, n. \quad (2.10)$$

## 2.2 Convergence Analysis

Let the exact solution of (2.4)-(2.5) $u$ satisfies following assumptions:

**Assumption 1:** $u \in C(\Omega) \cap C^4(\Omega \setminus \Gamma)$, and $\partial^4 u / \partial \theta^4$ is bounded in $\Omega$ with a constant $C_0$ such that $\sup\limits_{r \in (0,1)} (1-r)^{j-\sigma} |(\partial^j u / \partial r^j)(t, r, \theta)| \leq C_0$, $1 \leq j \leq 4$, where $\sigma \in (0, 2)$.

**Assumption 2:** There exists a constant $C_1$ such that $\omega(d) \equiv \sup\limits_{\text{dist}(P,Q) \leq d} |u(t, P) - u(t, Q)| \leq C_1 d^\sigma$ holds for any points $P$ and $Q$ near $\Gamma$.

**Assumption 3:** There exists a constant $C_2$ such that $|\partial^i u(t, r, \theta) / \partial t^i| \leq C_2$, $0 \leq i \leq 3$ for any $(t, r, \theta) \in \bar{D} = \Omega \times (0, T)$.

When the exact solution satisfies Assumptions 1-3 and exhibits singularity at the boundary of the domain, as $h$ approaches zero, the truncation error of Scheme (2.5)-(2.10) tends to infinity, thus the scheme is inconsistent. Therefore, traditional convergence results for non-singular problems cannot be applied.

For $P = (r_i, \theta_j) \in \Omega_h$, let $P_W = (r_{i-1}, \theta_j)$, $P_E = (r_{i+1}, \theta_j)$, $P_S = (r_i, \theta_{j-1})$, and $P_N = (r_i, \theta_{j+1})$. $P_W$, $P_E$, $P_S$, and $P_N$ are called the neighbors of $P$. Let $\Im < 1/4$ be a very small positive number. This paper arranges the grid points, $\Omega_h = \left( \bigcup\limits_{i=1}^{I} \Omega_h^{(i)} \right) \bigcup \Omega_h^{(0)}$, where $I = \lfloor \Im / h \rfloor$, which represents the largest integer not exceeding $\Im / h$,

$$\Omega_h^{(1)} = \{ P \in \Omega_h \mid \text{at least one neighbor of } P \in \Gamma_h \},$$

$$\Omega_h^{(i)} = \{ P \in \Omega_h \setminus \bigcup_{j=1}^{i-1} \Omega_h^{(j)} \mid \text{at least one neighbor of } P \in \Omega_h^{(i-1)} \}, \quad 2 \leq i \leq I,$$

and $\Omega_h^{(0)} = \Omega_h \setminus \bigcup\limits_{j=1}^{I} \Omega_h^{(j)}$, The number of points in $\Omega_h^{(i)}$ and $\Omega_h^{(0)}$ are $m_i$ and $m_0$ respectively ($m_0 + m_1 + \cdots + m_I = mn + 1$).

Rewrite (2.5) as

$$(1 + (1-a)\Delta t \wp_h) U_{i,j}^k = (1 - a\Delta t \wp_h) U_{i,j}^{k-1} + a\Delta t f_{i,j}^{k-1} + (1-a)\Delta t f_{i,j}^k, \quad k = 1, 2, 3, \ldots \quad (2.11)$$

Let $\mathbf{U}^k = (U_1^k, U_2^k, \ldots, U_N^k)^t$ be an unknown N-dimensional vector, $N = mn + 1$ according to the above grid arrangement order, representing the approximate value of the exact solution $u$ at time $k$. Then the linear system takes the form:

$$(I + (1-a)\Delta t \aleph) \mathbf{U}^k = (I - a\Delta t \aleph) \mathbf{U}^{k-1} + a\Delta t \mathbf{F}^{k-1} + (1-a)\Delta t \mathbf{F}^k, \quad k = 1, 2, 3, \ldots,$$

where $\aleph$ is an $N \times N$ matrix, and $\mathbf{F}^k$ is a known N-dimensional vector determined by the nonlinear term $f(t,r,\theta)$ and the boundary values $\gamma(t,\theta)$.

If for a matrix $A = (a_{i,j})$, any $i$, $j$ have $a_{i,j} \geq 0$, then use the notation $A \geq O$. Similarly, if $B - A \geq O$, then $B \geq A$.

We rearrange the grid points and denote:

$$\tilde{\mathbf{U}}^k = \left( \tilde{U}_{0,0}^k, \tilde{U}_{1,0}^k, \ldots, \tilde{U}_{m,0}^k, \tilde{U}_{1,1}^k, \ldots, \tilde{U}_{m,1}^k, \ldots, \tilde{U}_{1,n-1}^k, \ldots, \tilde{U}_{m,n-1}^k \right)^t,$$

as a new unknown $N$-dimensional vector. This is a new arrangement of $\mathbf{U}^k$. Here $U_{i,j}^k$ represents an approximate value of $u(k\Delta t, r_i, \theta_j)$. Then the linear system is written as:

$$(I + (1-a)\Delta t \tilde{\aleph}) \tilde{\mathbf{U}}^k = (I - a\Delta t \tilde{\aleph}) \tilde{\mathbf{U}}^{k-1} + a\Delta t \tilde{\mathbf{F}}^{k-1} + (1-a)\Delta t \tilde{\mathbf{F}}^k, \quad k = 1,2,3\ldots, \quad (2.12)$$

where $\tilde{\aleph}$ is the following $N \times N$ matrix:

$$\tilde{\aleph} = \begin{pmatrix} a_{0,0} & \alpha_0 & \alpha_1 & \alpha_2 & \cdots & \cdots & \alpha_{n-2} & \alpha_{n-1} \\ \beta_0 & A_0 & C_0 & 0 & \cdots & \cdots & 0 & B_0 \\ \beta_1 & B_1 & A_1 & C_1 & \ddots & & \ddots & 0 \\ \vdots & 0 & \ddots & \ddots & \ddots & \ddots & & \vdots \\ \vdots & \vdots & \ddots & \ddots & \ddots & \ddots & \ddots & \vdots \\ \vdots & \vdots & & \ddots & \ddots & \ddots & \ddots & 0 \\ \beta_{n-2} & 0 & \ddots & & \ddots & B_{n-2} & A_{n-2} & C_{n-2} \\ \beta_{n-1} & C_{n-1} & 0 & \cdots & \cdots & 0 & B_{n-1} & A_{n-1} \end{pmatrix}_{N \times N},$$

The elements of $\tilde{\aleph}$ are:

$$\alpha_j = \left( -\frac{4(1-a)}{nh_1^2}, 0, \cdots, 0 \right)_{1 \times m}, \quad j = 0,1,2,\cdots,n-1,$$

$$\beta_j = \left( -\frac{2(1-a)r_{\frac{1}{2}}}{r_1 h_1 (h_1 + h_2)}, 0, \cdots, 0 \right)_{m \times 1}^t, \quad j = 0,1,2,\cdots,n-1,$$

$\alpha_j$ and $\beta_j$ are $m$-dimensional vectors.

$$B_j = C_j = \begin{pmatrix} -\frac{1-a}{2r_1^2 \mu^2} & 0 & \cdots & 0 \\ 0 & -\frac{1-a}{2r_2^2 \mu^2} & \ddots & \vdots \\ \vdots & \ddots & \ddots & 0 \\ 0 & \cdots & 0 & -\frac{1-a}{2r_m^2 \mu^2} \end{pmatrix}_{m \times m}, \quad j = 0,1,2,\cdots,n-1,$$

$$A_j = \begin{pmatrix} a_1^{(j)} & c_1 & 0 & \cdots & \cdots & 0 \\ b_2 & a_2^{(j)} & c_2 & \ddots & & \vdots \\ 0 & \ddots & \ddots & \ddots & \ddots & \vdots \\ \vdots & \ddots & \ddots & \ddots & \ddots & 0 \\ \vdots & & \ddots & b_{m-1} & a_{m-1}^{(j)} & c_{m-1} \\ 0 & \cdots & \cdots & 0 & b_m & a_m^{(j)} \end{pmatrix}_{m \times m}, \quad j = 0,1,2,\cdots,n-1,$$

where

$$a_i^{(j)} = \frac{1}{\Delta t} + \frac{2(1-a)}{r_i(h_i + h_{i+1})}\left(\frac{r_{i+\frac{1}{2}}}{h_{i+1}} + \frac{r_{i-\frac{1}{2}}}{h_i}\right) + \frac{2(1-a)}{r_i^2 \mu^2}$$

$$+ 2(1-a)c_{i,j}, \quad i = 1,2,\cdots,m,$$

$$b_i = -\frac{2(1-a)r_{i-\frac{1}{2}}}{r_i h_i(h_i + h_{i+1})}, \quad i = 2,3,\cdots,m,$$

$$c_i = -\frac{2(1-a)r_{i+\frac{1}{2}}}{r_i h_{i+1}(h_i + h_{i+1})}, \quad i = 1,2,\cdots,m-1.$$

Since $\aleph$ is similar to $\tilde{\aleph}$ through a permutation matrix, it suffices to prove that $\tilde{\aleph}$ is an $M$-matrix. In fact, it can be easily verified that $\tilde{\aleph}$ is an irreducible matrix. Let $\tilde{\aleph} = (\tilde{w}_{i,j})$, and observe that

$$(\alpha_j)_i \leq 0, \quad (\beta_j)_i \leq 0, \quad b_i < 0, \quad c_i < 0, \quad a_{0,0} > 0, \quad a_i^{(j)} > 0$$

holds true, and note that $\tilde{w}_{i,i} \geq \sum_{j \neq i} |\tilde{w}_{i,j}|, \quad 1 \leq i \leq N$.

Let

$$E^k = \mathbf{u}^k - \mathbf{U}^k,$$

$$e^k = u(t+\Delta t, P) - u(t, P),$$

$$|\mathbf{e}^k| = \left(\max_{k(k\Delta t \leq T)} |e_1^k|, \ldots, \max_{k(k\Delta t \leq T)} |e_N^k|\right)^t,$$

$$|\mathbf{u}^k - \mathbf{U}^k| = \left(|u_1^k - U_1^k|, \ldots, |u_N^k - U_N^k|\right)^t.$$

**Lemma 2.1:** Suppose the exact solution $u$ of (2.1)-(2.3) satisfies Assumption 3, then there exists a positive constant $C$ independent of $h$, $\mu$, and $\Delta t$ such that

$$\sup_{K\Delta t \leq i} \|E^k\|_\infty \leq C(\tilde{\kappa}_1(a)\Delta t + \tilde{\kappa}_2(a)\Delta t^2).$$

**Proof:**

According to (2.11) and (2.12), the global error related to (2.12) can be written as

$$E^{k+1} = \mathbf{e}^{k+1} + RE^k,$$

where
$$R = \left(I + (1-a)\Delta t \wp L\right)^{-1}\left(I - a\Delta t \wp L\right).$$

The following recursive relationship is obtained:
$$E^{k+1} = \sum_{i=0}^{k} R^{k-i} \mathbf{e}^{k+1}.$$

Therefore, if the power of the transfer operator $R$ remains uniformly bounded, i.e.,
$$\| R^i \|_\infty \leq C, \quad i = 1, \cdots, n.$$

Then this lemma is immediately proven.□

The above argument proves the following conclusion.

**Theorem 2.1:** Suppose the exact solution $u$ of equations (2.1)-(2.3) satisfies Assumptions 1-3. Let $\mathbf{U}^k$ be the solution of scheme (2.5)-(2.10), and $\mathbf{u}^k$ be the vector of $u$ at grid points, where $t = k\Delta t$. Also assume that for some positive constant $M_0$, $\mu^2 \leq M_0 h$. Then if $\alpha = p\sigma < 2$, there exists a positive constant $c$ independent of $h$, $\mu$, and $\Delta t$ such that
$$\max\left|\mathbf{u}^k - \mathbf{U}^k\right| \leq c\left(\kappa_1(p)h^\alpha + \kappa_2(p)h^2 + \mu^2 + \tilde{\kappa}_1(a)\Delta t + \tilde{\kappa}_2(a)\Delta t^2\right),$$
for all $0 \leq k \leq K$.

Furthermore, if $p\sigma = 2$, then there also exists a positive constant $c$ independent of $h$, $\mu$, and $\Delta t$ such that
$$\max\left|\mathbf{u}^k - \mathbf{U}^k\right| \leq c\left(\kappa_1(p)h^2 |\log h| + \kappa_2(p)h^2 + \mu^2 + \tilde{\kappa}_1(a)\Delta t + \tilde{\kappa}_2(a)\Delta t^2\right),$$
for all $0 \leq k \leq K$.

The aforementioned $\kappa_1(p)$, $\kappa_2(p)$ are positive constants depending only on $p$ and increasing with respect to $p\sigma=2$. $\tilde{\kappa}_1(\alpha)$, $\kappa_2(p)$ are positive constants depending only on $\alpha$ and increasing with respect to $\alpha$.

It is noted that the proof of Theorem 2.1 can adopt a similar argument method as in Zhang et al.'s previous articles [46]. In fact, if the problem is static (i.e., not involving time changes), then this result is the main conclusion in Zhang et al.'s previous articles [46]. Therefore, the proof process is omitted here.

Theorem 2.1: Suppose the exact solution     of equations (2.1)-(2.3) satisfies Assumptions 1-3. Let $\mathbf{U}^k$ be the solution of scheme (2.5)-(2.10), and $\mathbf{u}^k$ be the vector of $u$ at grid points, where $t = k\Delta t$. Also assume that for some positive constant $M_0$, $\mu^2 \leq M_0 h$. Then if $\alpha = p\sigma < 2$, there exists a positive constant $c$ independent of $h$, $\mu$, and $\Delta t$ such that
$$\text{mean}\left|\mathbf{u}^k - \mathbf{U}^k\right| \leq c\left(\kappa_1(p)h^\alpha + \kappa_2(p)h^2 + \mu^2 + \tilde{\kappa}_1(a)\Delta t + \tilde{\kappa}_2(a)\Delta t^2\right),$$
for all $0 \leq k \leq K$.

Furthermore, if $p\sigma = 2$, then there also exists a positive constant $c$ independent of $h$, $\mu$, and $\Delta t$ such that

$$\text{mean}\left|\mathbf{u}^k - \mathbf{U}^k\right| \leq c\left(\kappa_1(p)h^2 |\log h| + \kappa_2(p)h^2 + \mu^2 + \tilde{\kappa}_1(a)\Delta t + \tilde{\kappa}_2(a)\Delta t^2\right),$$

for all $0 \leq k \leq K$.

The aforementioned $\kappa_1(p)$, $\kappa_2(p)$ are positive constants depending only on $p$ and increasing with respect to $p\sigma=2$. $\tilde{\kappa}_1(\alpha)$, $\kappa_2(p)$ are positive constants depending only on $\alpha$ and increasing with respect to $\alpha$.

It can be easily seen from Theorem 2.1 that Theorem 2.2 holds.

## 2.3 Numerical Experiment 1

In this section, the paper provides specific equations and initial conditions, along with some numerical results to discuss the effectiveness of the algorithm. The equation form and initial conditions are as follows:

$$\frac{\partial u}{\partial t} - \frac{1}{r}\frac{\partial}{\partial r} r \frac{\partial u}{\partial r} + \frac{1}{r^2}\frac{\partial^2 u}{\partial \theta^2} = f(t,r,\theta), \quad t \in (0,T),$$

$$u(t,\tilde{r},\theta) = e^{-t}, \quad t \in (0,T),$$

$$u(0,r,\theta) = r^2 + r^2(1-r)^{1/2},$$

$$f(t,r,\theta) = -e^{-t}r^2 - 4e^{-t} - 4(1-r)^{1/2} + \frac{5r}{2(1-r)^{1/2}} + \frac{r^2}{4(1-r)^{3/2}},$$

where $(r,\theta) \in \Omega = \{(r,\theta) \mid 0 < r < 1, 0 \leq \theta < 2\pi\}$ indicates that it is defined on a circular disk and the exact solution is $u(t,r,\theta) = e^{-t}r^2 + r^2(1-r)^{1/2}$ 。

For $m > 0$, let $h = 1/(m+1)$ and $n = [2m\pi]$, which is the largest integer less than $2m\pi$, and take $r_i = \varphi(i*h)$ and $\mu = 2\pi/n$. The stretching equation in Zhang et al.'s article is shown in formula (2.13) [43].

$$\varphi(s) = 1 - (1-s)^{\tilde{p}+1} \quad (0 \leq s \leq 1). \tag{2.13}$$

Firstly, using the grid division method from an article by Zhang et al. to find the optimal parameter in the stretching function (2.13) [43]. Figure 2.2 shows that when $a = 0.4$, $t = 0.1$ and, the optimal is approximately 0.6.

As shown in Table 2.1, when 当 , , , , , the error results of applying stretching function (2.4), applying stretching function (2.13) at  , and using a uniform grid are compared. The results indicate that constructing a non-uniform grid using stretching function (2.4) is much better than refining a uniform grid. The table also shows that the non-uniform grid division method used in this paper performs better than the grid point division determined by the optimal parameters using stretching function (2.13). Additionally, the self-adaptive method in this paper does not even need to use the best parameters to achieve the best performance (minimum maximum error) in the comparison, and it converges more stably with an accuracy improvement of about 7 times, which demonstrates the superiority of the algorithm in this paper. From the table, we can also see that when becomes smaller, the ratio in the table is roughly bounded by a constant, indicating that the numerical calculation results verify the theoretical results of error analysis.

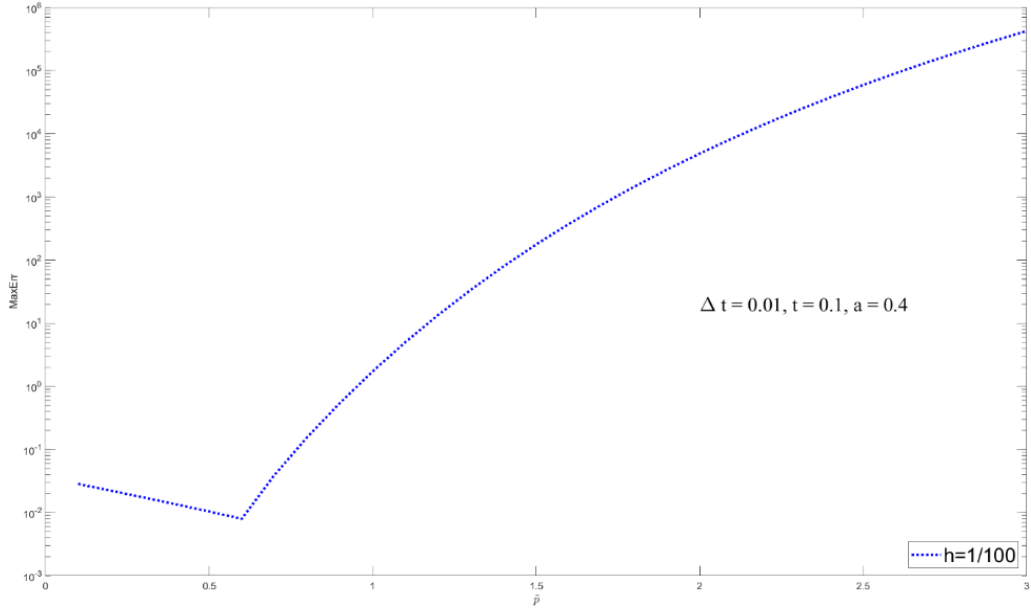

Figure 2.2 Maximum errors with respect to $\tilde{p}$ plot

Table 2.1 Numerical experiment 1 results comparison

| m | Maxerr(Stretching Function (2.4)) | $\dfrac{\text{Maxerr}}{h^{\alpha}+\Delta t}$ | Maxerr(Radius Non-Adaptive) | Maxerr(Stretching Function (2.13)) | Maxerr(Non-Adaptive) |
|---|---|---|---|---|---|
| 19 | 1.190276784e-02 | 2.310119305e-01 | 7.508413631e-02 | 2.733391989e-02 | 8.588423022e-02 |
| 39 | 6.353368652e-03 | 2.312356650e-01 | 5.418496890e-02 | 1.633457570e-02 | 9.046392948e-02 |
| 59 | 4.347949340e-03 | 2.117925643e-01 | 4.491190447e-02 | 1.194674879e-02 | 9.201692995e-02 |
| 79 | 3.311442609e-03 | 1.908677353e-01 | 3.913809524e-02 | 9.541207762e-03 | 9.279838657e-02 |
| 99 | 2.679448695e-03 | 1.721936019e-01 | 3.511282614e-02 | 8.005272818e-03 | 9.326884658e-02 |

When $m=99$ and $t=0.1$, the exact solution of the equation system is plotted in Figure 2.3, and when $m=99$, $t=0.1$, $\Delta t=0.01$, $a=0.4$, and $p=1$, the approximate solution of the equation system is plotted in Figure 2.4. Observing the three-dimensional graphs of the exact solution and the approximate solution, it can be seen that the shape of the approximate solution is very close to the exact solution, indicating good simulation results.

Figure 2.5 shows the results of how the maximum error changes with    and when  ,  , and  . Figure 2.6 shows the results of how the maximum error changes with and   when  ,  , and  . In both figures, the maximum error first decreases and then increases as    increases, and it also shows that the maximum error decreases as the grid division becomes finer. This indicates that the algorithm is convergent. Furthermore, there exists an optimal value of  , approximately 1.1, which can minimize the maximum error.

Figure 2.7 shows the results of how the maximum error changes with    and when , , and . The figure shows that the maximum error first increases and then decreases as aa increases, and there exists an optimal value of , approximately 0.5, which minimizes the maximum error when and . Additionally, the figure shows that the maximum error decreases as the grid division becomes finer, indicating that the algorithm is convergent.

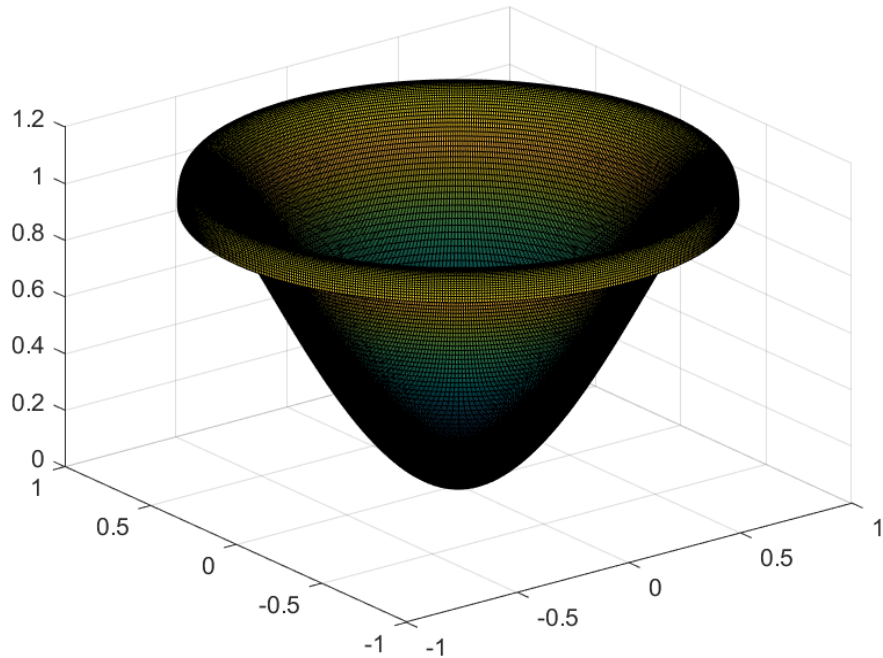

**Figure 2.3 Exact solution 3D plot**

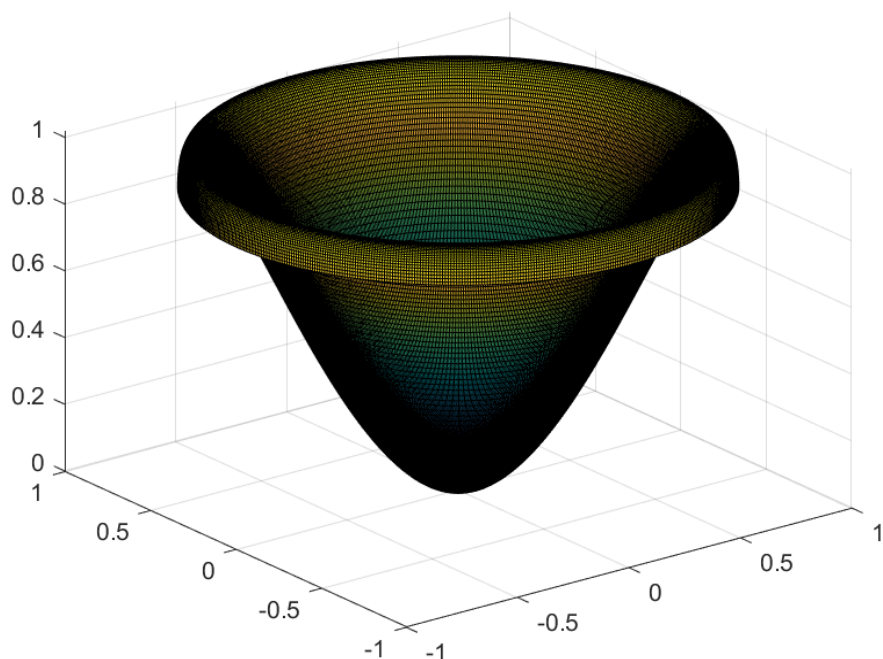

**Figure 2.4 Approximate solution 3D plot**

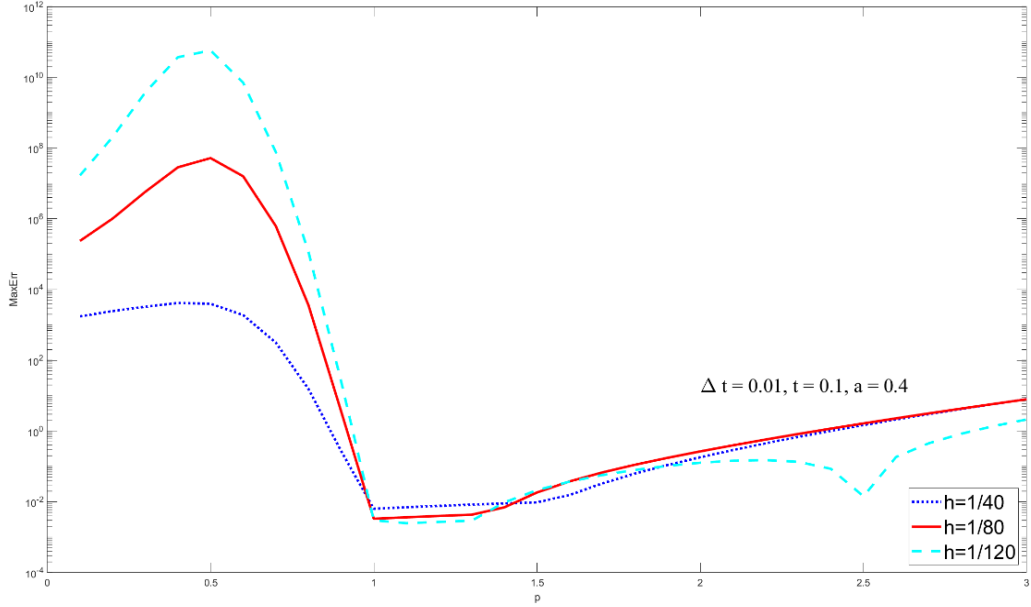

**Figure 2.5 Maximum errors with respect to $p$ and $h$ plot**

The error distribution of the approximate solution under the conditions $t = 0.1$, $h = 1/100$, $a = 0.4$, and $p = 1$ is shown in Figure 2.8, which demonstrates an extremely uniform error distribution with a very small overall numerical simulation error.

Figure 2.9 shows the results of how the maximum error varies with $h$ and $p$ when $t = 0.1$, $\Delta t = 0.02$, and $a = 0.4$. Figure 2.10 presents the results for the maximum error varying with $h$ and $p$ when $t = 0.1$, $\Delta t = 0.02$, and $a = 0.5$. Figure 2.11 illustrates the results of the maximum error varying with $h$ and $a$ when $t = 0.1$, $\Delta t = 0.02$, and $p = 1.1$. From these figures, it can be observed that the optimal value of $p$ consistently approaches 1.1; when the time step size is relatively large, the optimal value of $a$ tends towards 0.4, indicating its suitability for fast and accurate calculations; whereas when the time step size is smaller, the optimal value of $a$ leans towards 0.5, suggesting its appropriateness for high-precision calculations requiring high-performance computing support. Additionally, from the various figures in this section, it can be seen that the algorithm unconditionally converges as the time step length decreases.

In summary, the choice of parameters has a significant impact on the numerical simulation accuracy of this algorithm. In practical applications, selecting appropriate parameters based on specific parameter value experiments can greatly improve the accuracy of numerical simulations.

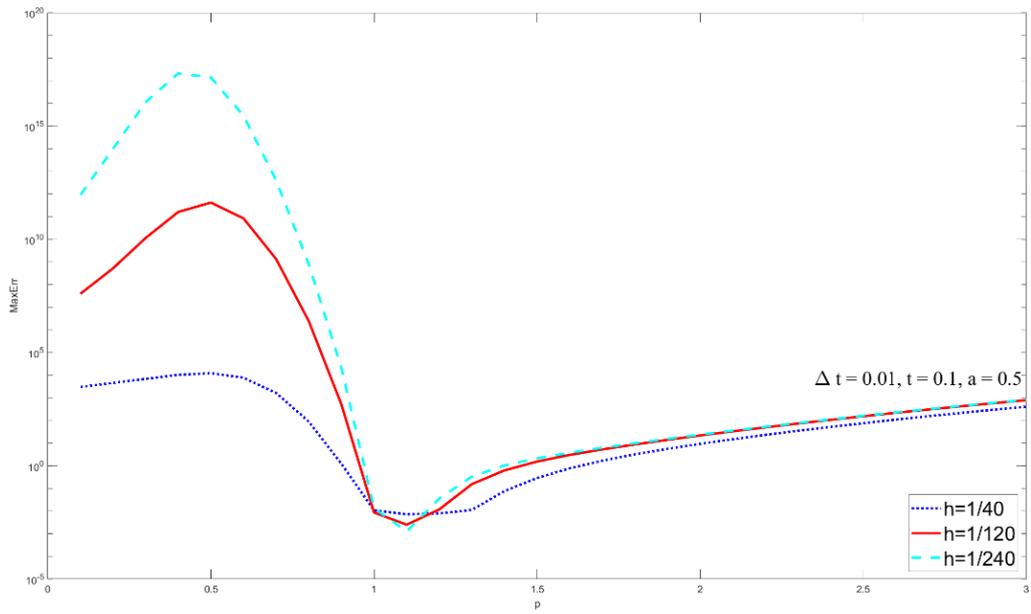

**Figure 2.6 Maximum errors with respect to $p$ and $h$ plot**

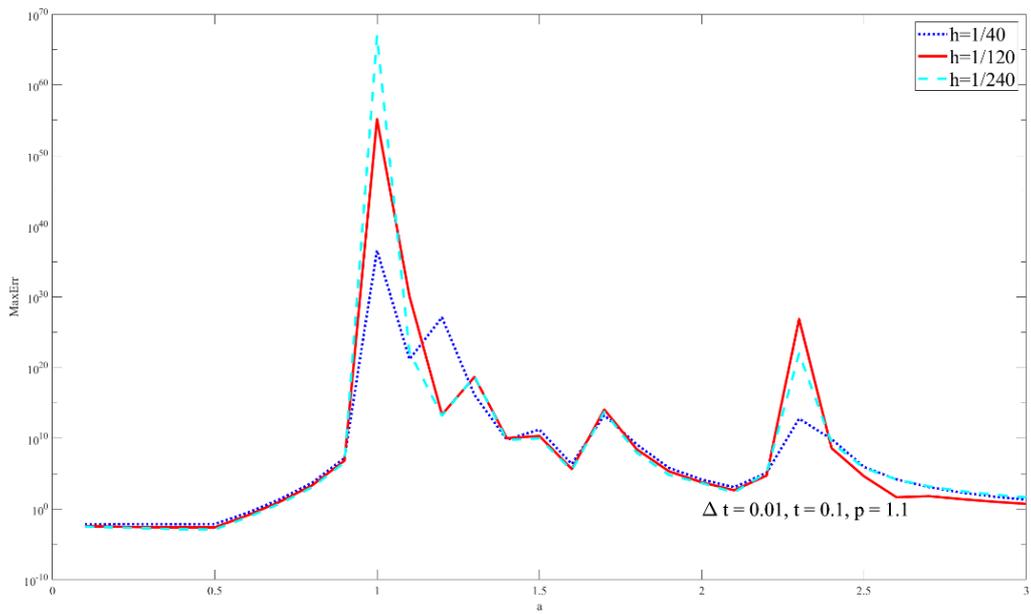

**Figure 2.7 Maximum errors with respect to $a$ and $h$ plot**

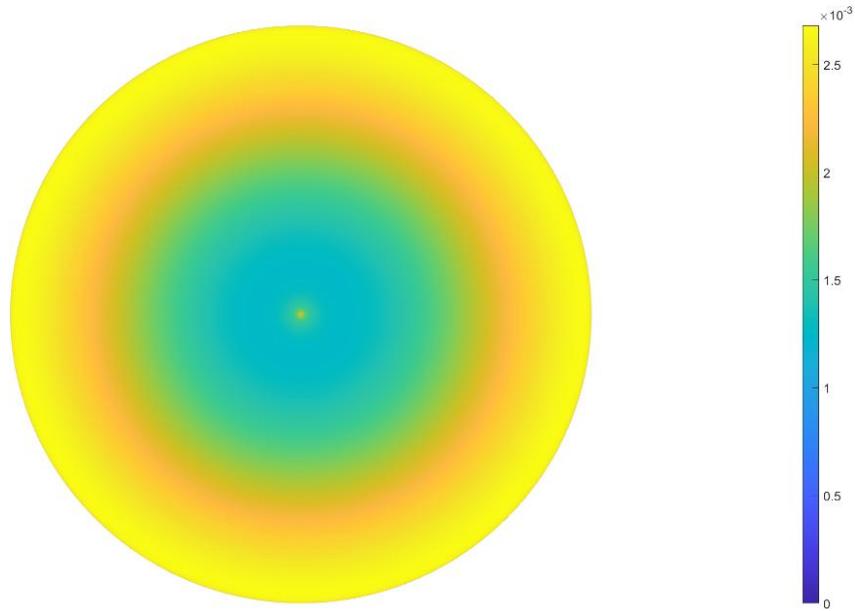

**Figure 2.8 Maximum errors distribution plot**

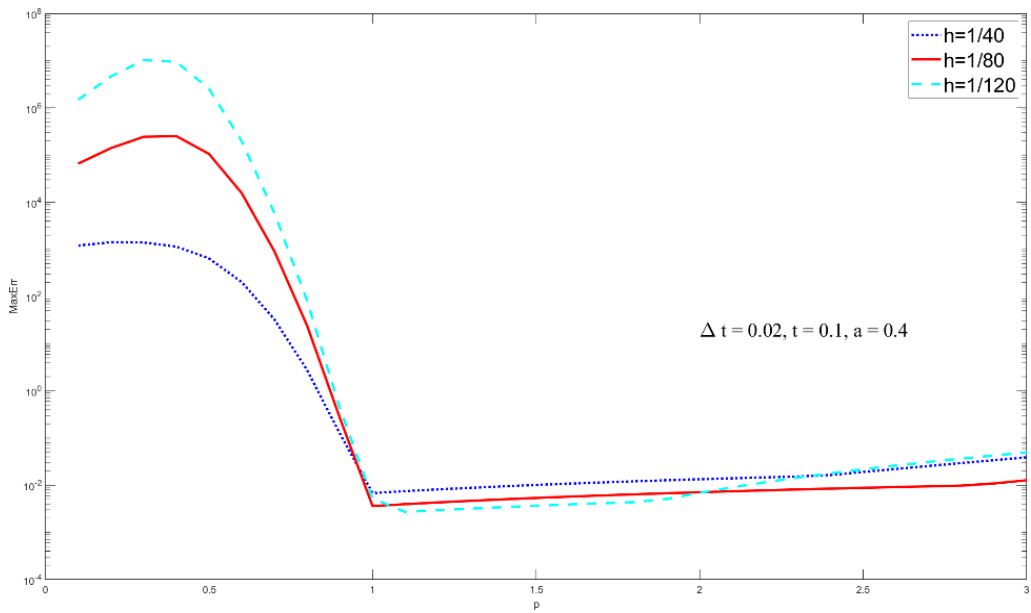

**Figure 2.9 Maximum errors with respect to $p$ and $h$ plot**

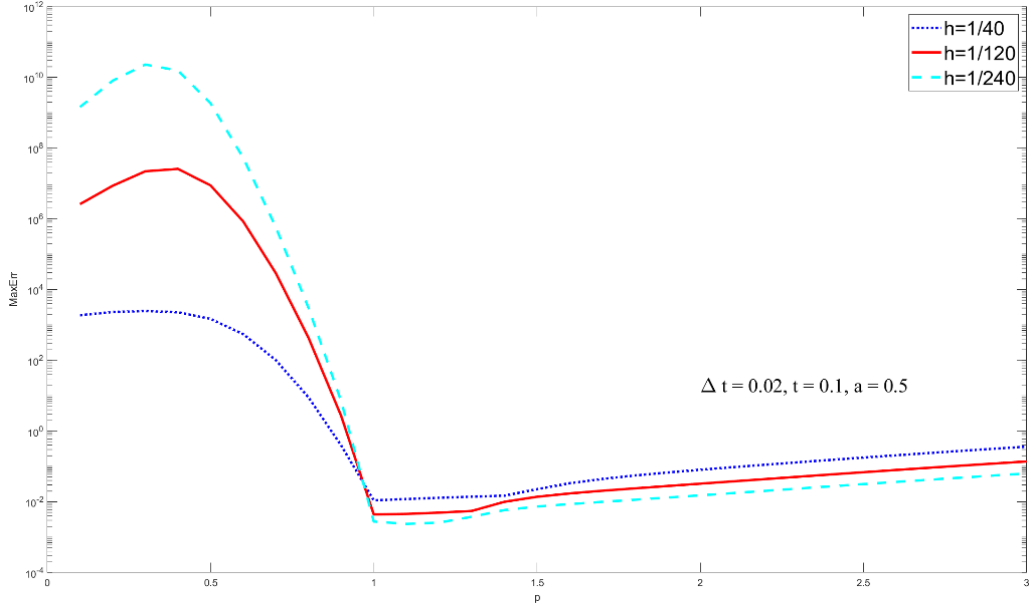

**Figure 2.10 Maximum errors with respect to** $p$ **and** $h$ **plot**

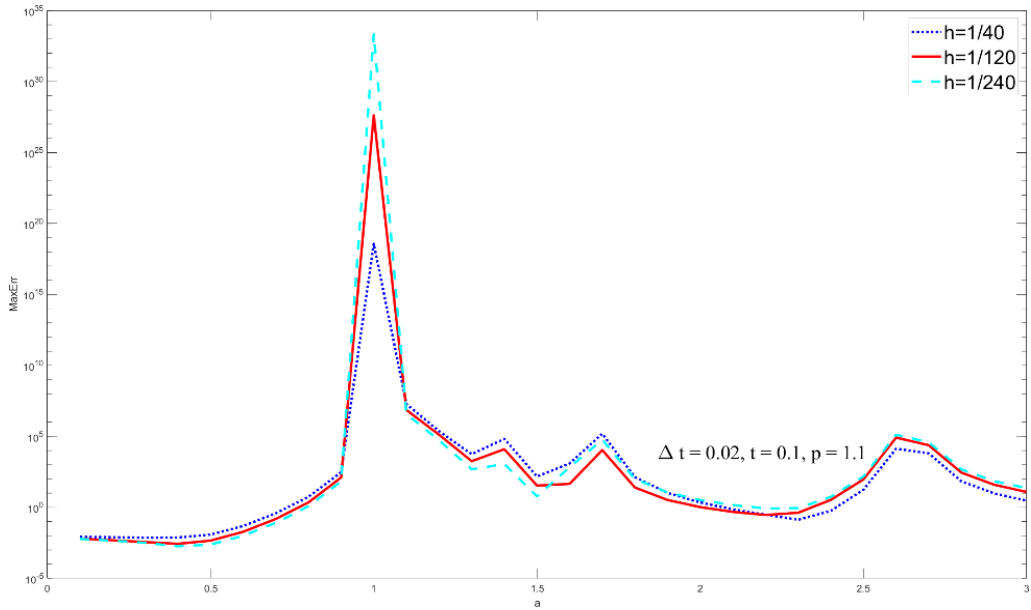

**Figure 2.11 Maximum errors with respect to** $a$ **and** $h$ **plot**

## 2.4 Numerical Experiment 2

In this section, a more complex form of the equation and its initial conditions are presented. In this case, the equation to be solved is related to angle, time, and radius, and involves more functional forms. Furthermore, this section provides some numerical results for analysis and discussion. The equation and initial conditions are as follows:

$$\frac{\partial u}{\partial t} - \left[\frac{1}{r}\frac{\partial}{\partial r}\left(r\frac{\partial u}{\partial r}\right) + \frac{1}{r^2}\frac{\partial^2 u}{\partial \theta^2}\right] + b(r,\theta)u = f(t,r,\theta), 0 < t < T,$$

$$u(t,\tilde{r},\theta) = \ln(1+e^t) + \sin\theta^2, \quad t \in (0,T),$$

$$b(r,\theta) = r^{1/2},$$

$$f(t,r,\theta) = \frac{r^2}{1+e^t}e^t - [4(1-r)^{1/2} - \frac{5}{2}r(1-r)^{-1/2}$$
$$-\frac{1}{4}r^2(1-r)^{-3/2} + 4\ln(1+e^t) + 4\sin\theta^2 + 2\cos 2\theta]$$
$$+r^{1/2}[r^2(1-r)^{1/2} + r^2\ln(1+e^{-t}) + \sin\theta^2 r^2],$$

where the expression $(r,\theta) \in \Omega = \{(r,\theta) | 0 < r < 1, 0 \le \theta < 2\pi\}$ indicates that it is defined on a disk and the exact solution is given by $u(t,r,\theta) = r^2(1-r)^{1/2} + r^2\ln(1+e^t) + \sin\theta^2 r^2$.

For $m > 0$, let $h = 1/(m+1)$ and $n = [2m\pi]$, which is the largest integer less than $2m\pi$. Define $r_i = \varphi(i*h)$ and $\mu = 2\pi/n$.

Firstly, we apply the mesh refinement method proposed in a paper by Zhang et al. [43] to find the optimal parameter $\tilde{p}$ for the stretching function (2.13). Figure 2.12 shows that when $a = 0.4$, $t = 0.00001$, $m = 99$, and $\Delta t = 0.000001$, the optimal $\tilde{p}$ is approximately 0.1.

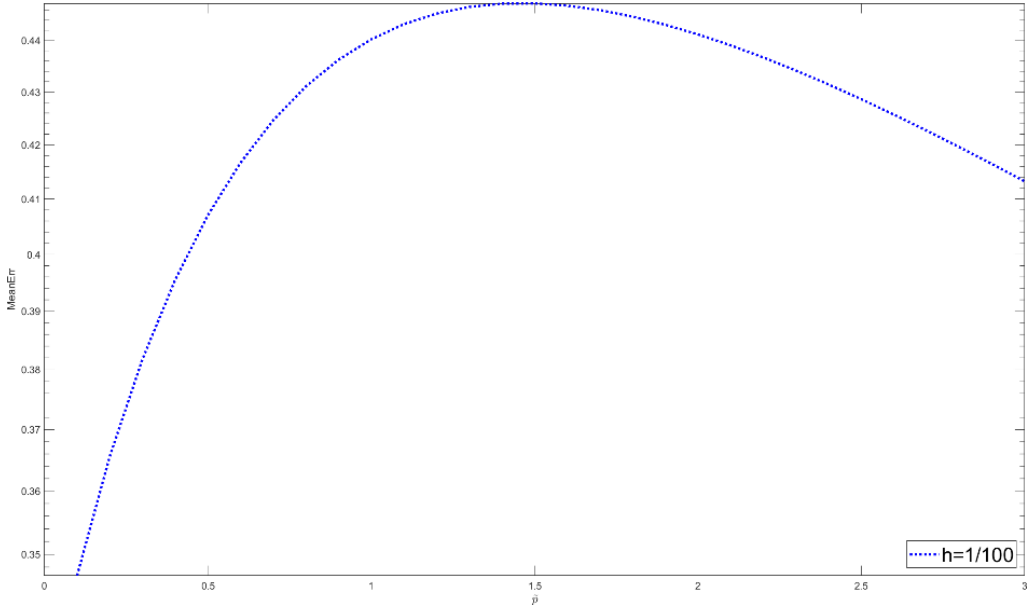

Figure 2.12 Mean errors with respect to $\tilde{p}$ plot

As shown in Table 2.2, when $t = 0.00001$, $\Delta t = 0.000001$, $p = 2.5$, $a = 0.4$, and $\alpha = p\sigma = 0.01$, the error results of using stretching function (2.4) are compared with those of using stretching function (2.13) with $\tilde{p}=0.1$. The results indicate that constructing non-uniform grids using stretching function (2.4) is significantly better than using uniform grid division. The table also shows that the grid division effect of the non-uniform grids applied in this paper is better than that determined by the optimal parameters for stretching function (2.13). Additionally, the self-adaptive method used in this paper achieves the best performance (minimum average error) without even needing to use the optimal parameters, and it converges more stably with an accuracy improvement of approximately one order of magnitude, demonstrating the superiority of the algorithm. From the table, we can also see that even as $h$ becomes smaller, the ratio of the values in the table remains roughly above a constant, which verifies the theoretical results of the error analysis.

The exact solutions of the equation set at $m = 99$ and $t = 0.00001$ are plotted alongside the approximate solutions of the equation set at $m = 99$, $t = 0.00001$, $\Delta t = 0.000001$, $a = 0.4$, and $p = 2.5$ in Figures 2.13 and 2.14. Observing the three-dimensional plots of the exact solution and the approximate solution, it can be seen that the shape of the approximate solution is very close to the exact solution, indicating good simulation results.

Table 2.2 Numerical experiment 2 results comparison

| m | Maxerr(Stretching Function (2.4)) | $\dfrac{\text{Maxerr}}{h^\alpha + \Delta t}$ | Maxerr(Radius Non-Adaptive) | Maxerr(Stretching Function (2.13)) | Maxerr(Non-Adaptive) |
|---|---|---|---|---|---|
| 19 | 2.427007318e-01 | 2.485503661e-01 | 3.229459692e-01 | 3.494925810e-01 | 4.185244520e-01 |
| 39 | 2.331405807e-01 | 2.404202005e-01 | 3.254346549e-01 | 3.496716459e-01 | 4.240253137e-01 |
| 59 | 2.310969078e-01 | 2.392774062e-01 | 3.254133754e-01 | 3.483576949e-01 | 4.258387535e-01 |
| 79 | 2.301393042e-01 | 2.389727872e-01 | 3.250357468e-01 | 3.473581122e-01 | 4.267418483e-01 |
| 99 | 2.297363332e-01 | 2.390865547e-01 | 3.246598033e-01 | 3.466898607e-01 | 4.272825320e-01 |

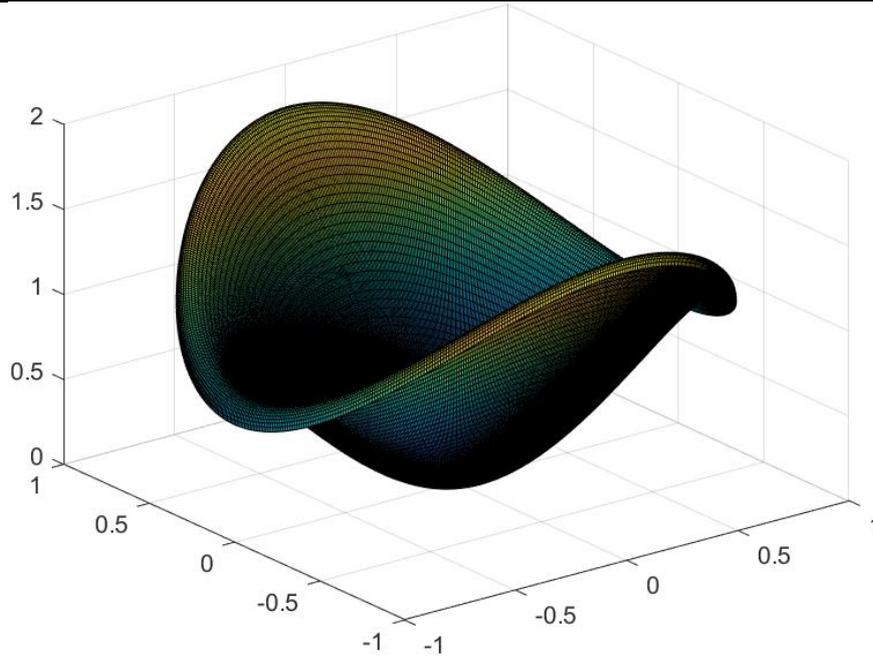

Figure 2.13 Exact solution 3D plot

Figure 2.15 shows the results of how the average error varies with $h$ and $p$ when $t = 0.00001$, $\Delta t = 0.000001$, and $a = 0.36$. The figure indicates that the average error first decreases and then increases as $p$ increases, and it also shows that the average error becomes smaller as the grid division becomes finer. This demonstrates that the algorithm is convergent. Furthermore, there exists an optimal value for $p$, approximately 0.2, which minimizes the average error.

Figure 2.16 displays the results of how the average error varies with $h$ and $a$ when $t = 0.00001$, $\Delta t = 0.000001$, and $p = 1$. The figure shows that the average error first decreases and then increases as $a$ increases, and there exists an optimal value for $a$ of approximately 0.36, which minimizes the average error when $t = 0.00001$ and $\Delta t = 0.000001$. Additionally, the figure shows that the average error becomes smaller as the grid division becomes finer, indicating that the algorithm is convergent. Combining the exploration of the optimal parameters for $a$ and $p$, it is found that in this relatively complex situation, a relatively smaller $a$ and $p$ can lead to smaller errors. This also indicates that the values of $a$ and $p$ have a significant impact on the numerical simulation results under different practical conditions. When these values are close to their optimal values, the accuracy of the numerical simulation is higher compared to non-adaptive algorithms. The complexity of the algorithm in this paper has not changed relative to non-adaptive algorithms, so selecting better parameters based on actual situations can greatly increase the accuracy of numerical simulations without changing the algorithm's efficiency.

The error distribution of the approximate solution under the conditions $t = 0.00001$, $\Delta t = 0.000001$, $h = 1/100$, $a = 0.36$, and $p = 0.2$ is shown in Figure 2.17. The error distribution is relatively uniform, with smaller middle errors, indicating good performance.

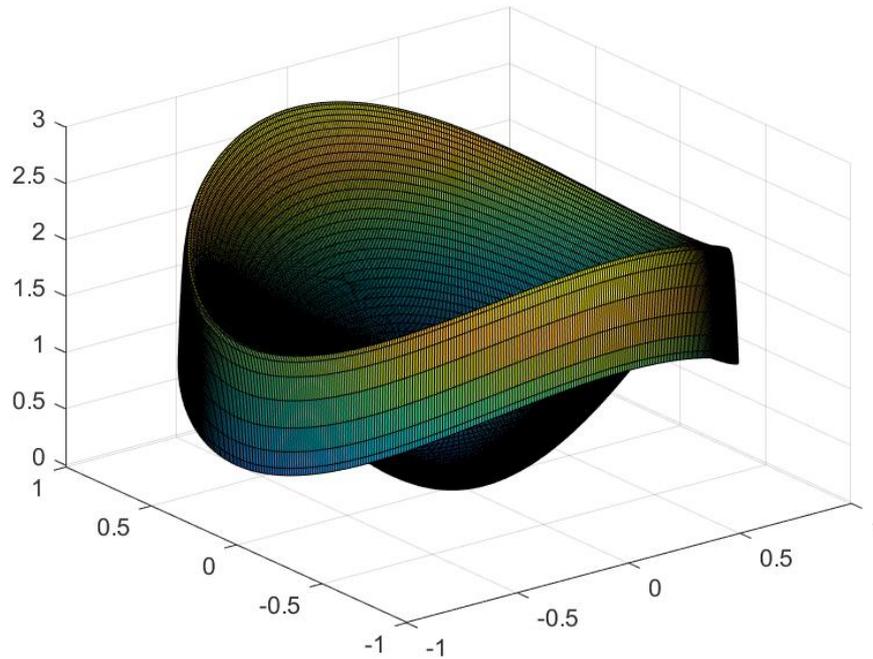

**Figure 2.14 Approximate solution 3D plot**

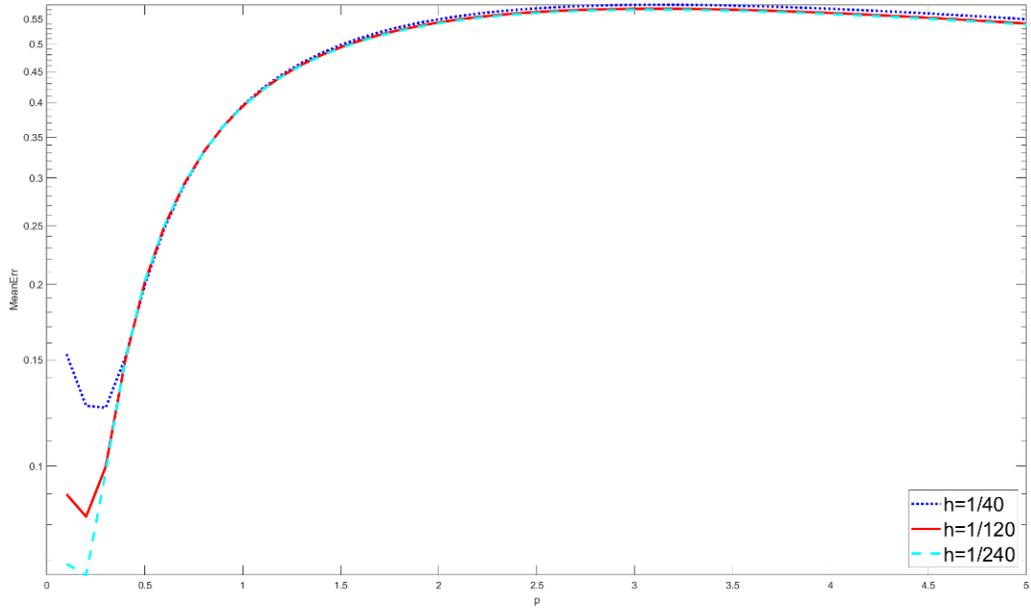

**Figure 2.15 Mean errors with respect to $p$ and $h$ plot**

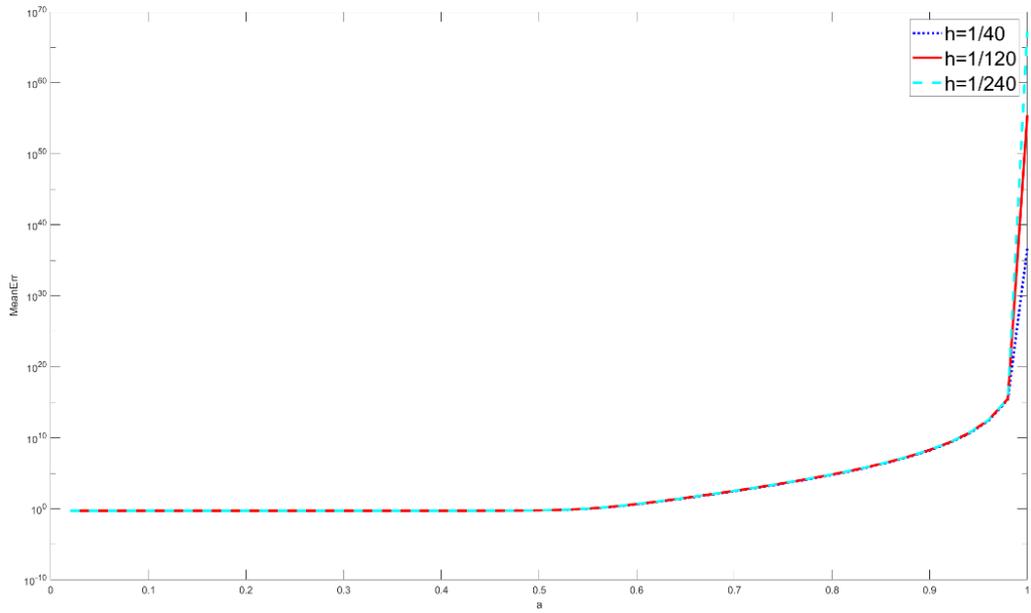

**Figure 2.16 Mean errors with respect to $a$ and $m$ plot**

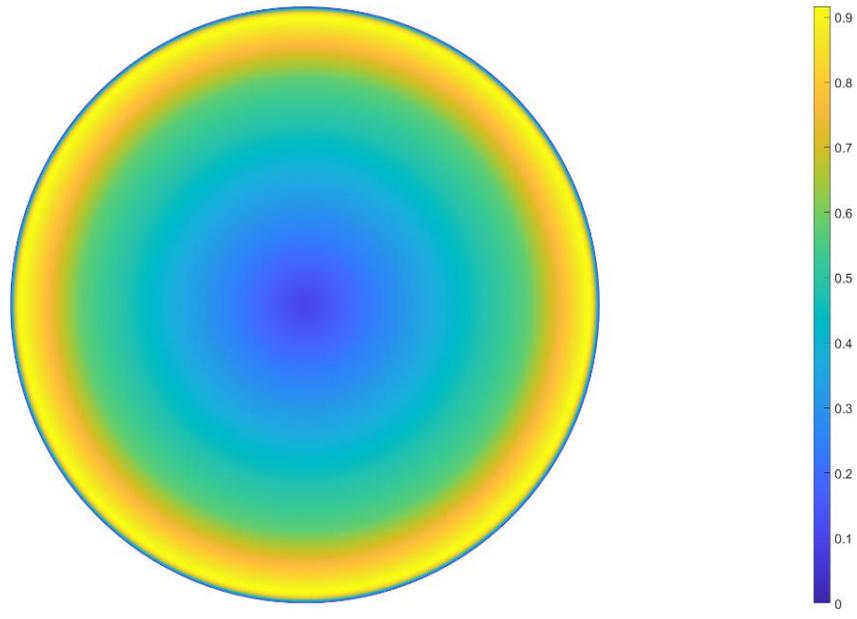

**Figure 2.17 Mean errors distribution plot**

# 3 Adaptive Numerical Simulation of Continuity Equations Derived from Boltzmann Equation

This chapter systematically explains the processing of the multi-dimensional adaptive numerical simulation framework in this paper for continuity equations, and a series of numerical experiments are conducted to verify the superiority of the algorithm. Additionally, certain analyses and discussions on the effects of the algorithm and the properties of parameters are provided.

## 3.1 Equation Form and Grid Division

The classical Boltzmann equation is given by:

$$\frac{\partial f}{\partial t} + \mathbf{v} \cdot \nabla_r f + \frac{\mathbf{F}}{m} \cdot \nabla_v f = I[f].$$

Through the quasi-neutral approximation, low-field condition assumption, and other simplifications, the continuity equations for electrons and holes (drift-diffusion model) can be derived as follows:

$$\frac{\partial \boldsymbol{n}}{\partial t} = \nabla \cdot (D_n \nabla \boldsymbol{n} - \mu_n \boldsymbol{n} \mathbf{E}) + G_n - R_n,$$

$$\frac{\partial \boldsymbol{p}}{\partial t} = \nabla \cdot (D_p \nabla \boldsymbol{p} - \mu_p \boldsymbol{p} \mathbf{E}) + G_p - R_p.$$

These equations are defined within the unit disk region in Cartesian coordinates: $\tilde{\Omega} = \{(x, y, t) \mid x^2 + y^2 < 1, 0 < t < T\}$. By performing a polar coordinate transformation on the original equations, we derive the continuity equations in polar coordinates on the unit disk:

$$\frac{\partial \boldsymbol{n}}{\partial t} = \left(\frac{1}{r}\frac{\partial}{\partial r}\left(rD_n \frac{\partial \boldsymbol{n}}{\partial r}\right) + \frac{1}{r^2}\frac{\partial}{\partial \theta}\left(D_n \frac{\partial \boldsymbol{n}}{\partial \theta}\right)\right)$$
$$-\left(\frac{1}{r}\frac{\partial}{\partial r}(r\boldsymbol{n}\mu_n \mathbf{E}) + \frac{1}{r}\frac{\partial}{\partial \theta}(\boldsymbol{n}\mu_n \mathbf{E})\right) + G_n - R_n,$$

$$\frac{\partial \boldsymbol{p}}{\partial t} = \left(\frac{1}{r}\frac{\partial}{\partial r}\left(rD_p \frac{\partial \boldsymbol{p}}{\partial r}\right) + \frac{1}{r^2}\frac{\partial}{\partial \theta}\left(D_p \frac{\partial \boldsymbol{p}}{\partial \theta}\right)\right)$$
$$-\left(\frac{1}{r}\frac{\partial}{\partial r}(r\boldsymbol{p}\mu_p \mathbf{E}) + \frac{1}{r}\frac{\partial}{\partial \theta}(\boldsymbol{p}\mu_p \mathbf{E})\right) + G_p - R_p.$$

Since $G_n$ and $R_n$ are typically directly provided by experimental data in physics, we denote $G_n - R_n$ as $GR(t, r, \theta)$, which is a known function. Given that the two equations have the same form, they can be written in the following general form for numerical solution:

$$\frac{\partial \boldsymbol{n}}{\partial t}(\frac{1}{r}\frac{\partial}{\partial r}\left(rD(r,\theta)\frac{\partial \boldsymbol{n}}{\partial r}\right)+\frac{1}{r^2}\frac{\partial}{\partial \theta}(D(r,\theta)\frac{\partial \boldsymbol{n}}{\partial \theta}))-(\frac{1}{r}\frac{\partial}{\partial r}(r\boldsymbol{n}\hat{\mu}(r,\theta)\mathbf{E}(r,\theta))$$
$$+\frac{1}{r}\frac{\partial}{\partial \theta}(\boldsymbol{n}\hat{\mu}(r,\theta)\mathbf{E}(r,\theta)))=GR(t,r,\theta), \quad 0<t<T, \quad (r,\theta)\in\Omega. \tag{3.1}$$

Then, the initial conditions and boundary conditions are given as follows:
$$\boldsymbol{n}(t,\tilde{r},\theta)=N(t,\theta), \quad 0<t<T, \quad (\tilde{r},\theta)\in\Gamma, \tag{3.2}$$
$$\boldsymbol{n}(0,r,\theta)=n_0(r,\theta), \quad (r,\theta)\in\Omega, \tag{3.3}$$
$$\hat{\mu}(r,\theta)\mathbf{E}(r,\theta)=E(r,\theta). \tag{3.4}$$

The domain of definition for the system of equations is: $\Omega=\{(r,\theta)\,|\,0<r<1, 0\leq\theta<2\pi\}$, $\Gamma=\{(\tilde{r},\theta)\,|\,\tilde{r}=1, 0\leq\theta<2\pi\}$.

The radial stretching equation is:
$$\upsilon(s)=\sin(\frac{s^p\pi}{2}), \quad 0\leq s\leq 1, \quad p>0.$$

The angular stretching equation is:
$$v(\iota)=2\pi\sin(\frac{(\frac{\iota}{2\pi})^q\pi}{2}), \quad 0\leq\iota\leq 2\pi, \quad q>0.$$

The time stretching equation is:
$$V(t)=T\sin(\frac{(\frac{t}{T})^l\pi}{2}), \quad 0\leq t\leq T, \quad l>0.$$

When $p>0$, the radial stretching equation satisfies $\upsilon(0)=0$, $\upsilon(1)=1$. When $q>0$, the angular stretching equation satisfies $v(0)=0$ and $v(2\pi)=2\pi$. When $l>0$, the time stretching equation satisfies $V(0)=0$ and $V(T)=T$.

Subsequently, the region is divided as follows:
$$h=\frac{1}{m+1}, \quad r_i=\upsilon(ih), \quad r_{i+\frac{1}{2}}=\frac{r_i+r_{i+1}}{2}, \quad i=0,1,2,\cdots,m+1,$$
$$h_i=r_i-r_{i-1}, \quad i=1,2,\cdots,m+1,$$
$$\tilde{\mu}=\frac{2\pi}{n}, \quad \theta_j=v(j\tilde{\mu}), \quad j=0,1,2,\cdots,n,$$
$$\mu_j=\theta_j-\theta_{j-1}, \quad j=1,2,\cdots,n, \quad \mu_0=\mu_n,$$
$$\mu_{j+\frac{1}{2}}=\frac{\mu_j+\mu_{j+1}}{2}, \quad j=0,1,2,\cdots,n,$$
$$\Delta\tilde{t}=\frac{1}{k}, \quad k=0,1,2,\cdots,K, \quad K\in Z^+,$$
$$t_k=V(k\Delta\tilde{t}), \quad k=0,1,2,\cdots,K,$$
$$\Delta t_{k+1}=t_{k+1}-t_k, \quad k=0,1,2,\cdots,K-1.$$

The regions after division in the time and spatial dimensions are shown in Figures 3.1 and 3.2, respectively.

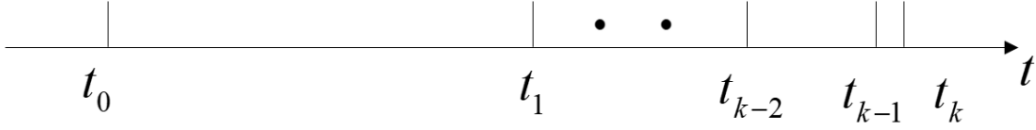

Figure 3.1 Plot of temporal lattice division

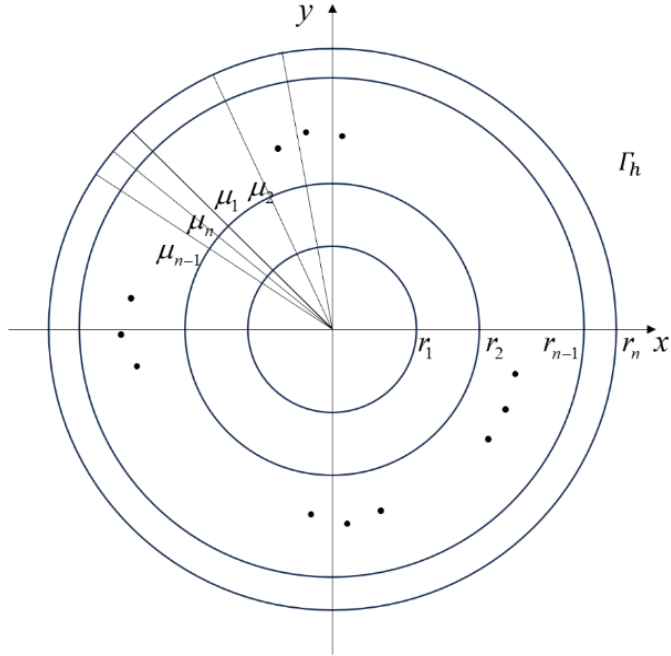

Figure 3.2 Plot of spatial lattice division

## 3.2 Discrete Format and Iterative Linear Equation System Derivation
### 3.2.1 Derivation of Non-Origin Point Discrete Format

The original equation is discretized using the adjustable coefficient adaptive Crank-Nicolson method in the time dimension direction and at non-origin points in the spatial dimension, as follows:

$$\frac{n_{i,j}^{k+1} - n_{i,j}^{k}}{\Delta t} + a\varpi_h n_{i,j}^{k} + (1-a)\varpi_h n_{i,j}^{k+1} = aGR_{i,j}^{k} + (1-a)GR_{i,j}^{k+1}, \quad k \in Z^{+}, \quad 0 < a < 1,$$

where $n_{i,j}^{k}$ is the approximate value of the exact solution $n(t_k, r_i, \theta_j)$ for equations (3.1)-(3.4) at the time layer $t_k = V(k\Delta \tilde{t})$, and $GR_{i,j}^{k} = f(t_k, r_i, \theta_j)$,

$$\varpi_h n_{i,j} = \begin{pmatrix} \dfrac{2}{r_i(h_i+h_{i+1})}\left(\dfrac{r_{i+\frac{1}{2}}D_{i+\frac{1}{2},j}}{h_{i+1}} + \dfrac{r_{i-\frac{1}{2}}D_{i-\frac{1}{2},j}}{h_i}\right) \\ + \dfrac{2}{r_i^2(\mu_j+\mu_{j+1})}\left(\dfrac{D_{i,j+\frac{1}{2}}}{\mu_{j+1}} + \dfrac{D_{i,j-\frac{1}{2}}}{\mu_j}\right) + E_{i,j}\left(\dfrac{1}{r_i\mu_{j+1}} + \dfrac{1}{h_{i+1}}\right) \end{pmatrix} n_{i,j}$$

$$-\dfrac{2r_{i-\frac{1}{2}}D_{i-\frac{1}{2},j}}{r_i h_i(h_i+h_{i+1})}n_{i-1,j} - \left(\dfrac{2r_{i+\frac{1}{2}}D_{i+\frac{1}{2},j}}{r_i h_{i+1}(h_i+h_{i+1})} + E_{i,j}\dfrac{1}{h_{i+1}}\right)n_{i+1,j} - \dfrac{2D_{i,j-\frac{1}{2}}}{r_i^2 \mu_j(\mu_j+\mu_{j+1})}n_{i,j-1}$$

$$-\left(\dfrac{2D_{i,j+\frac{1}{2}}}{r_i^2 \mu_{j+1}(\mu_j+\mu_{j+1})} + E_{i,j}\dfrac{1}{r_i\mu_{j+1}}\right)n_{i,j+1}, \quad i=1,2,\cdots,m, \quad j=0,1,2,\cdots,n-1.$$

### 3.2.2 Derivation of Discrete Format at the Origin

Assume that:

1. Axial Symmetry: Near the origin, the physical quantity $n$ is symmetric about $\theta$, i.e., $\dfrac{\partial n}{\partial \theta}=0$.

2. Coefficient Continuity: $D(r,\theta)$, $\mu_n(r,\theta)$, and $E(r,\theta)$ are bounded as $r \to 0$. Set $D(0,\theta)=D_0$ and $,E(0,\theta)=E_0$ to define the values of these functions at $r=0$.

Under the above assumptions, we derive the discrete format at the origin using an integral method. We average the values around the origin as $n_1 = \dfrac{1}{n}\sum_{j=0}^{n-1} n_{1,j}$ and integrate the original equation over a small annular region $r \in [0, \Delta r/2]$ and $\theta \in [0, 2\pi]$:

$$\int_0^{2\pi}\int_0^{\Delta r/2} \dfrac{\partial n}{\partial t} r\, dr\, d\theta = \text{other}.$$

Time Term Integration:

$$\int_0^{2\pi}\int_0^{\Delta r/2} \dfrac{\partial n}{\partial t} r\, dr\, d\theta \approx \pi\left(\dfrac{\Delta r}{2}\right)^2 \dfrac{n_0^{k+1}-n_0^k}{\Delta t},$$

where $n_0$ is the value at the origin, and $k$ is the time step.

Diffusion Term Integration:

$$\int_0^{2\pi}\int_0^{\Delta r/2} \dfrac{1}{r}\dfrac{\partial}{\partial r}\left(rD(r,\theta)\dfrac{\partial n}{\partial r}\right) r\, dr\, d\theta = 2\pi r D(r,\theta)\dfrac{\partial n}{\partial r}\bigg|_{r=\Delta r/2}.$$

At $r = \Delta r/2$:

Using central difference approximation for the derivative: $\dfrac{\partial n}{\partial r}\bigg|_{r=\Delta r/2} \approx \dfrac{n_1 - n_0}{\Delta r}$.

Assume $D(r,\theta)$ is approximately equal to $D_0$ (the origin value) within the interval $r \in [0, \Delta r/2]$.

Therefore, the diffusion term integral result is:

$$2\pi\left(\frac{\Delta r}{2}\right)D_0\frac{\boldsymbol{n}_1-\boldsymbol{n}_0}{\Delta r}=\pi\Delta r D_0\frac{\boldsymbol{n}_1-\boldsymbol{n}_0}{\Delta r}.$$

Drift Term Integral:

$$-\int_0^{2\pi}\int_0^{\Delta r/2}\frac{1}{r}\frac{\partial}{\partial r}(r n\mu_n(r,\theta)E(r,\theta))r dr d\theta=-2\pi r n\mu_n(r,\theta)E(r,\theta)|_{r=\Delta r/2}.$$

At $r=\Delta r/2$, take the average value $n_{1/2}=\dfrac{n_0+n_1}{2}$, and assume $\mu_n$ and $E$ are $\mu_{n_0}$ and $E_0$. Therefore, the drift term integral result is:

$$-2\pi\left(\frac{\Delta r}{2}\right)\left(\frac{\boldsymbol{n}_0+\boldsymbol{n}_1}{2}\right)\mu_{n_0}E_0=-\pi\Delta r\mu_{n_0}E_0\left(\frac{\boldsymbol{n}_0+\boldsymbol{n}_1}{2}\right).$$

Angular Derivative Term Integral:

By the symmetry assumption $\dfrac{\partial n}{\partial \theta}=0$, the angular derivative term integral is zero.

Source Term Integral:

$$\int_0^{2\pi}\int_0^{\Delta r/2} GR(t,r,\theta)r dr d\theta \approx \pi\left(\frac{\Delta r}{2}\right)^2 GR(t,0).$$

Combine all integral terms and divide by the control volume $\pi(\Delta r/2)^2$ to obtain the final form of the discrete format:

$$\frac{\boldsymbol{n}_0^{k+1}-\boldsymbol{n}_0^k}{\Delta t}=\frac{4D_0}{\Delta r^2}(\boldsymbol{n}_1^{k+1}-\boldsymbol{n}_0^{k+1})-\frac{2\mu_{n_0}E_0}{\Delta r}\left(\frac{\boldsymbol{n}_0^{k+1}+\boldsymbol{n}_1^{k+1}}{2}\right)+GR(t,0).$$

It can be seen that:

$$\varpi_h \boldsymbol{n}_{0,0}=\left(\frac{4D_{0,0}}{h_1^2}+\frac{E_{0,0}}{h_1}\right)\boldsymbol{n}_{0,0}-\left(\frac{4D_{0,0}}{n h_1^2}-\frac{E_{0,0}}{n h_1}\right)\sum_{j=0}^{n-1}\boldsymbol{n}_{1,j},$$

where $\boldsymbol{n}_{i,j}^k$ is the approximate value of the exact solution $\boldsymbol{n}(t_k,r_i,\theta_j)$ for equations (3.1)-(3.4) at the time layer $t_k=V(k\Delta\tilde{t})$ and $GR_{i,j}^k=f(t_k,r_i,\theta_j)$.

The equations at the origin and non-origin points satisfy the following initial and boundary conditions:

$$\boldsymbol{n}_{0,j}^k=\boldsymbol{n}_{0,0}^k,\quad \boldsymbol{n}_{m+1,j}^k=N_j^k,\quad j=0,1,2,\cdots,n,$$

$$\boldsymbol{n}_{i,n}^k=\boldsymbol{n}_{i,0}^k,\quad \boldsymbol{n}_{i,-1}^k=\boldsymbol{n}_{i,n-1}^k,\quad i=0,1,2,\cdots,m+1,$$

$$\boldsymbol{n}_{i,j}^0=\boldsymbol{n}_0(r_i,\theta_j),\quad i=0,1,2,\cdots,m+1;\ j=0,1,2,\cdots,n.$$

### 3.2.3 Iterative Linear Equation System Derivation

Establish a linear equation system to solve as follows:

$$\begin{pmatrix} a_{0,0} & \alpha_0 & \alpha_1 & \alpha_2 & \cdots & \cdots & \alpha_{n-2} & \alpha_{n-1} \\ \beta_0 & A_0 & C_0 & 0 & \cdots & \cdots & 0 & B_0 \\ \beta_1 & B_1 & A_1 & C_1 & \ddots & & & 0 \\ \vdots & 0 & \ddots & \ddots & \ddots & \ddots & & \vdots \\ \vdots & \vdots & \ddots & \ddots & \ddots & \ddots & \ddots & \vdots \\ \vdots & \vdots & & \ddots & \ddots & \ddots & \ddots & 0 \\ \beta_{n-2} & 0 & & \ddots & & B_{n-2} & A_{n-2} & C_{n-2} \\ \beta_{n-1} & C_{n-1} & 0 & \cdots & 0 & 0 & B_{n-1} & A_{n-1} \end{pmatrix} \begin{pmatrix} \boldsymbol{n}_{0,0}^{k+1} \\ \boldsymbol{n}_{1,0}^{k+1} \\ \boldsymbol{n}_{2,0}^{k+1} \\ \vdots \\ \boldsymbol{n}_{m,0}^{k+1} \\ \boldsymbol{n}_{1,1}^{k+1} \\ \vdots \\ \boldsymbol{n}_{m,1}^{k+1} \\ \vdots \\ \boldsymbol{n}_{1,n-1}^{k+1} \\ \vdots \\ \boldsymbol{n}_{m,n-1}^{k+1} \end{pmatrix} = \begin{pmatrix} e_{0,0}^k \\ e_{1,0}^k \\ e_{2,0}^k \\ \vdots \\ e_{m,0}^k \\ e_{1,1}^k \\ \vdots \\ e_{m,1}^k \\ \vdots \\ e_{1,n-1}^k \\ \vdots \\ e_{m,n-1}^k \end{pmatrix},$$

$$a_{0,0} = \frac{1}{\Delta t_{k+1}} + (\frac{4D_{0,0}}{h_1^2} + \frac{E_{0,0}}{h_1})(1-a),$$

$$\alpha_j = \left( -(1-a)(\frac{4D_{0,0}}{nh_1^2} - \frac{E_{0,0}}{nh_1}), 0, \cdots, 0 \right)_{1 \times m},$$

$$\beta_j = \left( -\frac{2(1-a)r_1 D_{\frac{1}{2},j}}{r_1 h_1 (h_1+h_2)}, 0, \cdots, 0 \right)_{m \times 1}^t, \quad j=0,1,2,\cdots,n-1,$$

$$B_j = \begin{pmatrix} -\frac{2(1-a)D_{1,j-\frac{1}{2}}}{r_1^2 \mu_j (\mu_j + \mu_{j+1})} & 0 & \cdots & 0 \\ 0 & -\frac{2(1-a)D_{2,j-\frac{1}{2}}}{r_2^2 \mu_j (\mu_j + \mu_{j+1})} & \ddots & \vdots \\ \vdots & \ddots & \ddots & 0 \\ 0 & \cdots & 0 & -\frac{2(1-a)D_{m,j-\frac{1}{2}}}{r_m^2 \mu_j (\mu_j + \mu_{j+1})} \end{pmatrix}_{m \times m},$$

$$j=0,1,2,\cdots,n-1,$$

$$C_j = \begin{pmatrix} -(1-a)(\dfrac{2D_{1,j+\frac{1}{2}}}{r_1^2 \mu_{j+1}(\mu_j+\mu_{j+1})} \\ +E_{1,j}\dfrac{1}{r_1 \mu_{j+1}}) & 0 & \cdots & 0 \\ 0 & -(1-a)(\dfrac{2D_{2,j+\frac{1}{2}}}{r_2^2 \mu_{j+1}(\mu_j+\mu_{j+1})} \\ +E_{2,j}\dfrac{1}{r_2 \mu_{j+1}}) & \ddots & \vdots \\ \vdots & \ddots & \ddots & 0 \\ 0 & \cdots & 0 & -(1-a)(\dfrac{2D_{m,j+\frac{1}{2}}}{r_m^2 \mu_{j+1}(\mu_j+\mu_{m+1})} \\ +E_{i,j}\dfrac{1}{r_i \mu_{j+1}}) \end{pmatrix}_{m \times m},$$

$j = 0,1,2,\cdots,n-1,$

$$A_j = \begin{pmatrix} a_1^{(j)} & c_1^{(j)} & 0 & \cdots & \cdots & 0 \\ b_2^{(j)} & a_2^{(j)} & c_2^{(j)} & \ddots & & \vdots \\ 0 & \ddots & \ddots & \ddots & \ddots & \vdots \\ \vdots & \ddots & \ddots & \ddots & \ddots & 0 \\ \vdots & & \ddots & b_{m-1}^{(j)} & a_{m-1}^{(j)} & c_{m-1}^{(j)} \\ 0 & \cdots & \cdots & 0 & b_m^{(j)} & a_m^{(j)} \end{pmatrix}_{m \times m}, \quad j=0,1,2,\cdots,n-1,$$

$$a_i^{(j)} = \frac{1}{\Delta t_{k+1}} + (1-a)\left( \frac{2}{r_i(h_i+h_{i+1})}\left( \frac{r_{i+\frac{1}{2}} D_{i+\frac{1}{2},j}}{h_{i+1}} + \frac{r_{i-\frac{1}{2}} D_{i-\frac{1}{2},j}}{h_i} \right) + \frac{2}{r_i^2(\mu_j+\mu_{j+1})}\left( \frac{D_{i,j+\frac{1}{2}}}{\mu_{j+1}} + \frac{D_{i,j-\frac{1}{2}}}{\mu_j} \right) + E_{i,j}\left( \frac{1}{r_i \mu_{j+1}} + \frac{1}{h_{i+1}} \right) \right),$$

$i=1,2,\cdots,m,$

$$e_{0,0}^k = \left( \frac{1}{\Delta t_{k+1}} - a(\frac{4D_{0,0}}{h_1^2} + \frac{E_{0,0}}{h_1}) \right) \mathbf{n}_{0,0}^k + a(\frac{4D_{0,0}}{nh_1^2} - \frac{E_{0,0}}{nh_1})\sum_{j=0}^{n-1} \mathbf{n}_{1,j}^k$$
$$+(1-a)GR_{0,0}^{k+1} + aGR_{0,0}^k,$$

$$c_i^{(j)} = -(1-a)(\frac{2r_{i+\frac{1}{2}} D_{i+\frac{1}{2},j}}{r_i h_{i+1}(h_i+h_{i+1})} + E_{i,j}\frac{1}{h_{i+1}}), \quad i=1,2,\cdots,m-1,$$

$$b_i^{(j)} = -\frac{2(1-a)r_{i-\frac{1}{2}} D_{i-\frac{1}{2},j}}{r_i h_i(h_i+h_{i+1})}, \quad i=2,3,\cdots,m,$$

$e_{i,j}^k = s_{i,j}^k, \quad i=1,2,\cdots,m-1, \quad j=0,1,2,\cdots,n-1,n,$

$$e_{m,j}^k = s_{m,j}^k + \frac{2ar_{m+\frac{1}{2}}D_{m+\frac{1}{2},j}}{r_m h_{m+1}(h_m + h_{m+1})} N_j^{k+1}, \quad j = 0,1,2,\cdots,n-1,n,$$

$$s_{i,j}^k = \left(\frac{1}{\Delta t_{k+1}} - a\left(\frac{2}{r_i(h_i + h_{i+1})}\left(\frac{r_{i+\frac{1}{2}}D_{i+\frac{1}{2},j}}{h_{i+1}} + \frac{r_{i-\frac{1}{2}}D_{i-\frac{1}{2},j}}{h_i}\right) + \frac{2}{r_i^2(\mu_j + \mu_{j+1})}\left(\frac{D_{i,j+\frac{1}{2}}}{\mu_{j+1}} + \frac{D_{i,j-\frac{1}{2}}}{\mu_j}\right) + E_{i,j}\left(\frac{1}{r_i\mu_{j+1}} + \frac{1}{h_{i+1}}\right)\right)\right)n_{i,j}^k$$

$$+ \frac{2ar_{i-\frac{1}{2}}D_{i-\frac{1}{2},j}}{r_i h_i(h_i + h_{i+1})} n_{i-1,j}^k + a\left(\frac{2r_{i+\frac{1}{2}}D_{i+\frac{1}{2},j}}{r_i h_{i+1}(h_i + h_{i+1})} + E_{i,j}\frac{1}{h_{i+1}}\right)n_{i+1,j}^k + \frac{2aD_{i,j-\frac{1}{2}}}{r_i^2 \mu_j(\mu_j + \mu_{j+1})} n_{i,j-1}^k$$

$$+ a\left(\frac{2D_{i,j+\frac{1}{2}}}{r_i^2 \mu_{j+1}(\mu_j + \mu_{j+1})} + E_{i,j}\frac{1}{r_i\mu_{j+1}}\right)n_{i,j+1}^k + (1-a)GR_{i,j}^{k+1} + aGR_{i,j}^k,$$

$$i = 1,2,\cdots,m, \quad j = 0,1,2,\cdots,n-1.$$

### 3.3 Numerical Experiment 1

In this section, the paper provides an example of a continuity equation along with initial conditions. In this scenario, the equation to be solved includes numerous functional forms and is relatively complex. Additionally, numerical experiments were conducted in this part, and analyses and discussions were carried out based on the results of these numerical experiments. Consider the following equation form and initial conditions:

$$\frac{\partial \boldsymbol{n}}{\partial t} - \left(\frac{1}{r}\frac{\partial}{\partial r}\left(rD(r,\theta)\frac{\partial \boldsymbol{n}}{\partial r}\right) + \frac{1}{r^2}\frac{\partial}{\partial \theta}\left(D(r,\theta)\frac{\partial \boldsymbol{n}}{\partial \theta}\right)\right) - \left(\frac{1}{r}\frac{\partial}{\partial r}(r\boldsymbol{n}\hat{\mu}(r,\theta)\mathbf{E}(r,\theta))\right.$$

$$+ \frac{1}{r}\frac{\partial}{\partial \theta}(\boldsymbol{n}\hat{\mu}(r,\theta)\mathbf{E}(r,\theta))) = GR(t,r,\theta), \quad 0 < t < T, \quad (r,\theta) \in \Omega,$$

$$\boldsymbol{n}(t,\tilde{r},\theta) = e^{-t} + \sin\theta^2, \quad 0 < t < T, \quad (\tilde{r},\theta) \in \Gamma,$$

$$\boldsymbol{n}(0,r,\theta) = r + r(1-r)^{1/2} + \sin\theta^2 r, \quad (r,\theta) \in \Omega,$$

$$\hat{\mu}(r,\theta)\mathbf{E}(r,\theta) = E(r,\theta) = r,$$

$$D(r,\theta) = r,$$

$$GR(t,r,\theta) = -e^{-t}r - (2e^{-t} + 2(1-r)^{1/2} - 2r(1-r)^{-1/2}$$

$$+ \frac{1}{4}r^2(1-r)^{-3/2} + 2\sin\theta^2 + 2\cos(2\theta)) + 3e^{-t}r + 3r(1-r)^{1/2}$$

$$- \frac{1}{2}r^2(1-r)^{-1/2} + 3r\sin\theta^2 + \sin 2\theta r,$$

where $(r,\theta) \in \Omega = \{(r,\theta) \mid 0 < r < 1, 0 \le \theta < 2\pi\}$ indicates that it is defined on a disk and has an exact solution $\boldsymbol{n}(t,r,\theta) = e^{-t}r + r(1-r)^{1/2} + \sin\theta^2 r$.

For $m > 0$, let $h = 1/(m+1)$ and $n = [2m\pi]$, which is the largest integer less than $2m\pi$, Take $r_i = \varphi(i*h)$ and $\tilde{\mu} = 2\pi/n$.

As shown in Table 3.1, when $t = 0.1$, $\Delta\tilde{t} = 0.01$, $a = 0.5$. $l = 1.5$, $p = 0.1$, $q = 1.9$, the maximum error results are compared using three adaptive stretching functions and non-adaptive forms for grid division. In this relatively simple example, the introduction of angle and difference non-adaptive methods does not significantly improve numerical simulation accuracy. However, the adaptive method applied in the radial direction shows a more noticeable improvement in solution precision. The numerical simulation under time non-adaptivity is scattered and unstable; thus, overall, the non-adaptive method in time is inferior to the adaptive method. This indicates that introducing an adaptive approach can enhance algorithm accuracy and stability, with an approximate 10% increase in precision.

**Table 3.1 Numerical experiment 1 results comparison**

| m | Maxerr(Uniform Adaptation) | Maxerr(Angle Non-adaptation) | Maxerr(Radius Non-adaptation) | Maxerr(Time Non-adaptation) | Maxerr(Difference Non-adaptation) |
|---|---|---|---|---|---|
| 19 | 3.905372902e-01 | 3.902743094e-01 | 2.928963823e-01 | 1.613165667e-01 | 3.905372902e-01 |
| 39 | 3.212570476e-01 | 3.211306935e-01 | 2.925879106e-01 | 1.614519888e-01 | 3.212570476e-01 |
| 59 | 2.893545221e-01 | 2.892995641e-01 | 2.918929966e-01 | 1.617956177e-01 | 2.893545221e-01 |
| 79 | 2.798637667e-01 | 2.798240126e-01 | 2.913269055e-01 | 1.620581914e-01 | 2.798637667e-01 |
| 99 | 2.755451493e-01 | 2.755166272e-01 | 2.909112100e-01 | 1.622574996e-01 | 2.755451493e-01 |

When $m = 99$, $t = 0.1$, the exact solutions of the equation system are plotted against the approximate solutions when $m = 99$, $t = 0.1$, $\Delta\tilde{t} = 0.01$, $a = 0.5$, $l = 1.5$, $p = 0.1$, as shown in Figures 3.3 and 3.4. Observing the three-dimensional plots of the exact and approximate solutions, it is evident that the shape of the approximate solution closely matches the exact solution, indicating good simulation results.

Figure 3.5 illustrates the average error variation with $h$ and $p$ when $t = 0.1$, $\Delta\tilde{t} = 0.01$, $a = 0.5$, $q = 1$, $l = 1$. The figure shows that the average error first increases and then decreases as $p$ increases, demonstrating that the average error decreases as the grid becomes finer. This indicates that the algorithm is convergent. There exists an optimal value of $p$, approximately 0.1, which minimizes the average error.

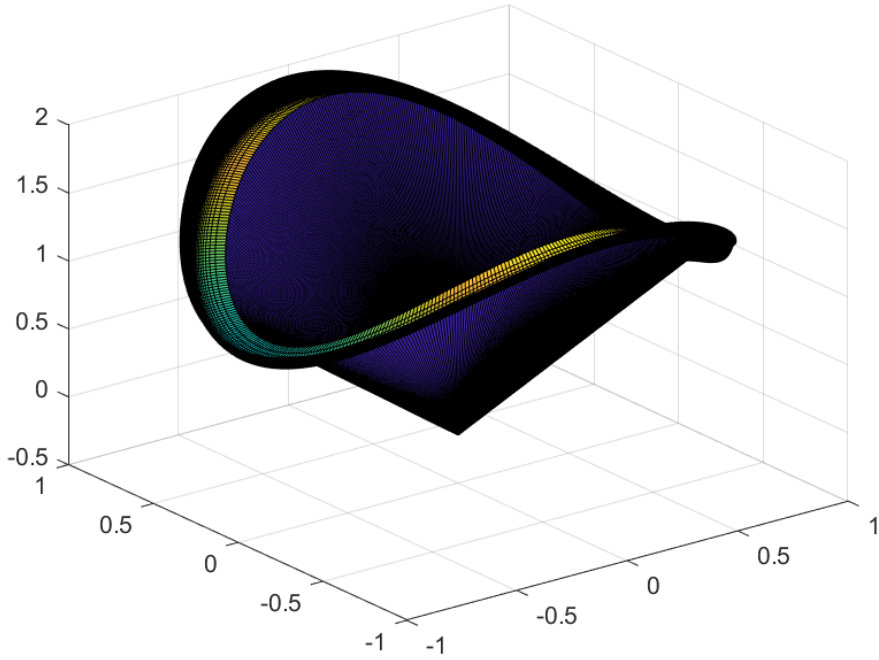
**Figure 3.3 Exact solution 3D plot**

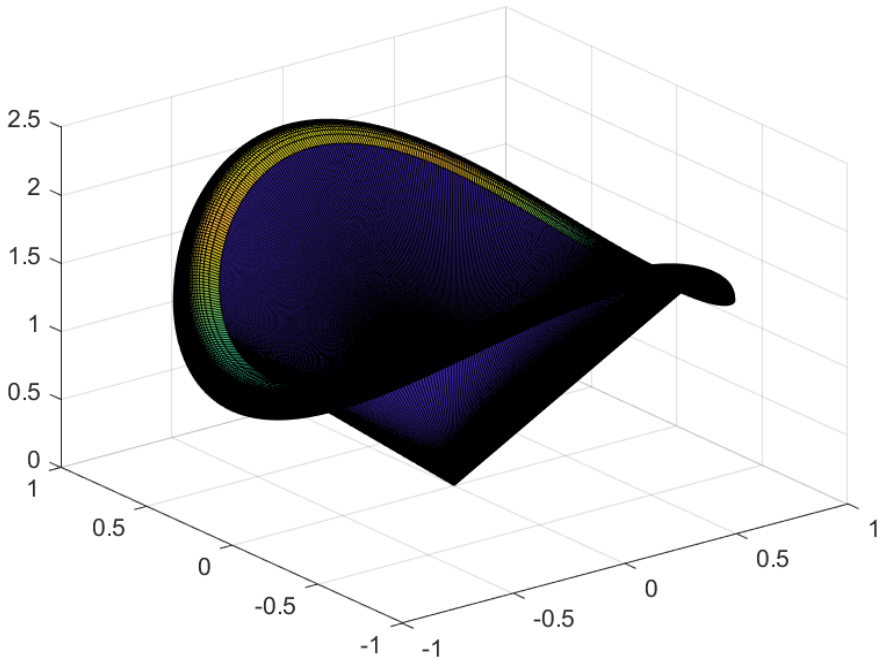
**Figure 3.4 Approximate solution 3D plot**

Figure 3.6 shows the results of how the average error changes with $h$ and $q$ when $t = 0.1$, $\Delta \tilde{t} = 0.01$, $a = 0.5$, $p = 0.1$, $l = 1$. The figure indicates that the average error first decreases and then increases as $q$ increases, and it also shows that the average error decreases as the grid becomes finer. This demonstrates that the algorithm is convergent. Furthermore, there exists an optimal value of $q$, approximately 1.9, which minimizes the average error.

Figure 3.7 displays the results of how the average error changes with $h$ and $l$ when $t = 0.1, \Delta \tilde{t} = 0.01, a = 0.5, p = 0.1, q = 1.9$. The figure shows that the average error first decreases and then increases as $l$ increases, and it also indicates that the average error decreases as the grid becomes finer. This demonstrates that the algorithm is convergent. Additionally, there exists an optimal value of $l$, approximately 1.5, which minimizes the average error.

Figure 3.8 shows the results of how the average error changes with $h$ and $a$ when $t = 0.1$, $\Delta \tilde{t} = 0.01$, $l = 1.5$, $p = 0.1$, $q = 1.9$. The figure indicates that the average error first increases and then decreases as $a$ increases, and it also shows that there is an optimal value of $a$ around 0.5, which minimizes the average error of the numerical solution when $t = 0.1$, $\Delta \tilde{t} = 0.01$. Additionally, the figure demonstrates that the average error decreases as the grid becomes finer, indicating that the algorithm is convergent.

By comprehensively analyzing the exploration of the optimal parameters for $a$ and $p$, it is found that in this relatively simple case, setting $a$ close to 0.5, $l$ close to 1.5, $q$ close to 1.9, and $p$ to a relatively small value can result in smaller errors. This also indicates that under different real-world conditions, the values of $a$ and $p$ significantly affect the results of numerical simulations. When these values are close to their optimal values, the accuracy of numerical simulations improves considerably compared to non-adaptive algorithms. The complexity of the algorithm presented in this paper does not change relative to non-adaptive algorithms; therefore, selecting better parameters based on actual conditions can greatly increase the accuracy of numerical simulations without changing the algorithm's efficiency.

At $t = 0.1$ and $\Delta \tilde{t} = 0.01$, the distribution of errors in the approximate solution under the conditions $h = 1/100$, $a = 0.5$, $l = 1.5$, $p = 0.1$, $q = 1.9$ is shown in Figure 3.9. The figure shows that the error distribution is relatively uniform. The angle stretching function and radial stretching function can significantly reduce the error at specific positions, thereby reducing the overall average error and improving the overall simulation effect.

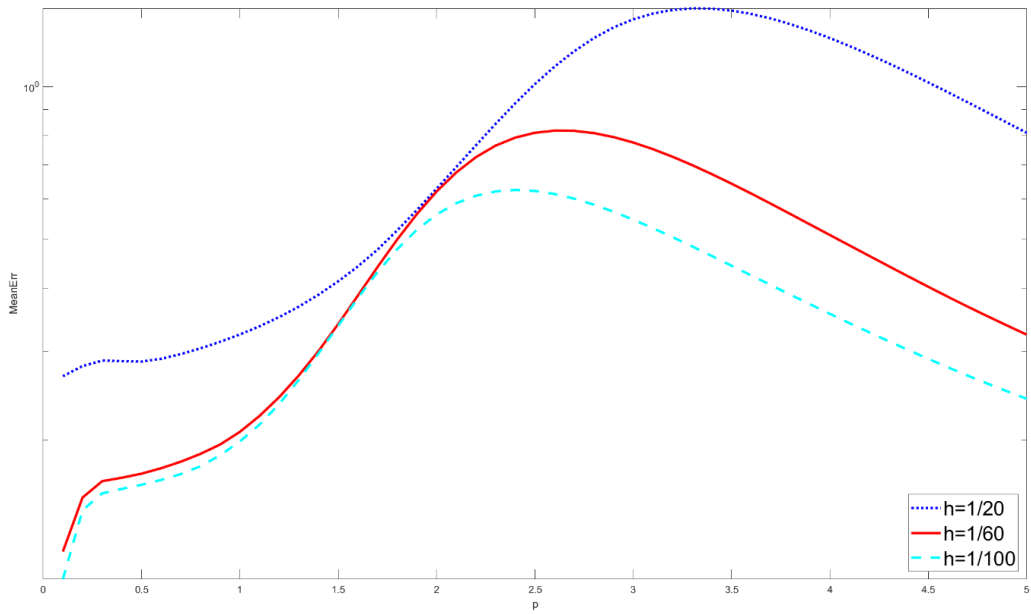

**Figure 3.5 Mean errors with respect to $p$ and $h$ plot**

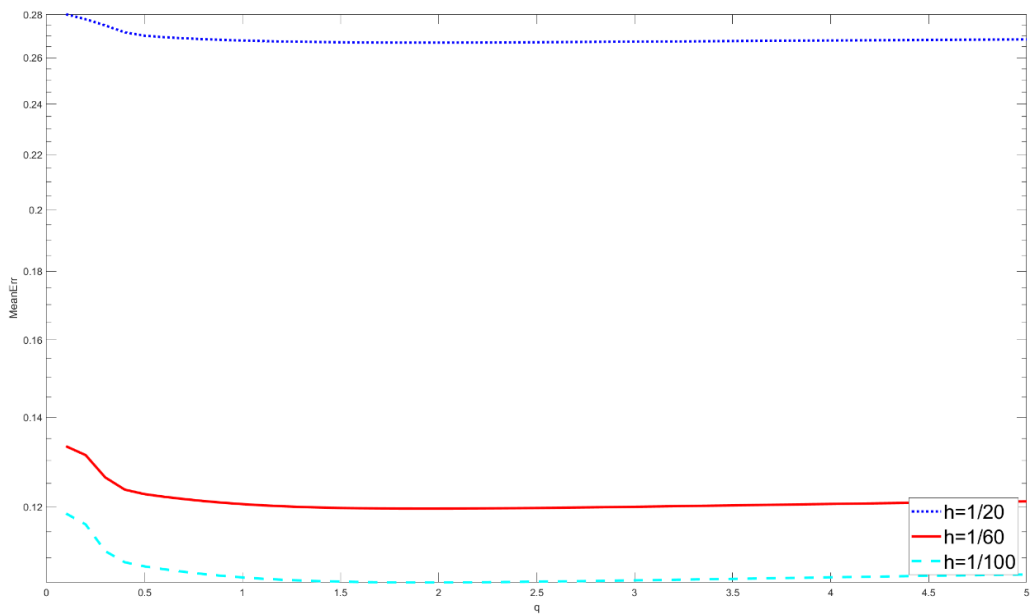

**Figure 3.6 Mean errors with respect to $q$ and $h$ plot**

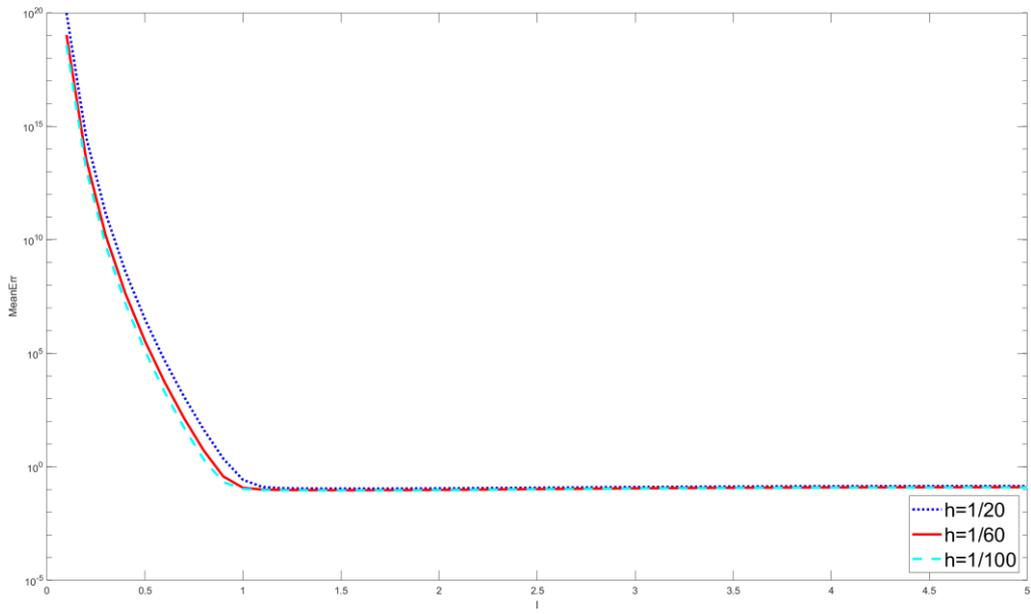

**Figure 3.7 Mean errors with respect to** $l$ **and** $h$ **plot**

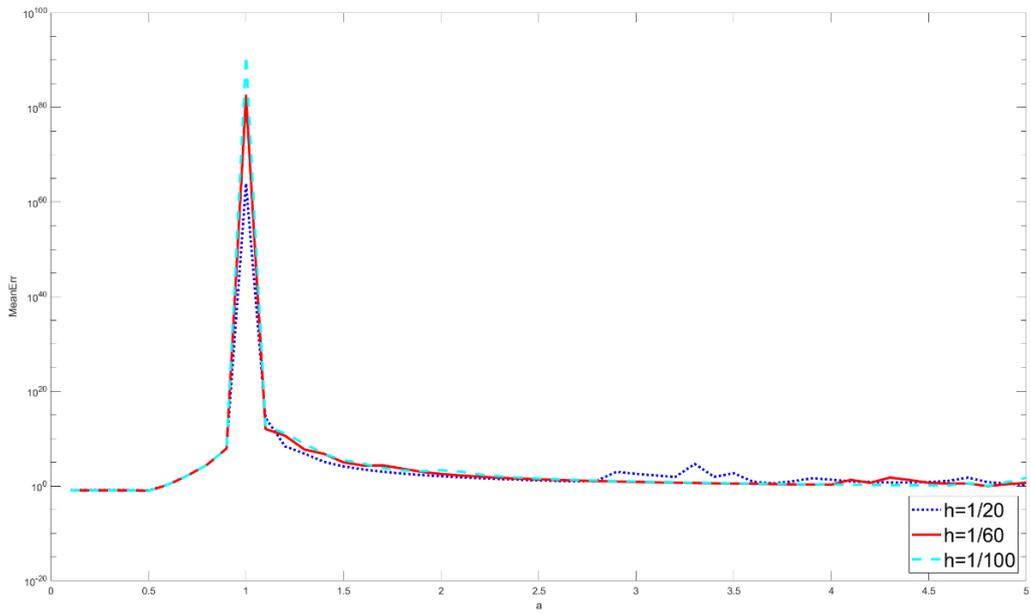

**Figure 3.8 Mean errors with respect to** $a$ **and** $h$ **plot**

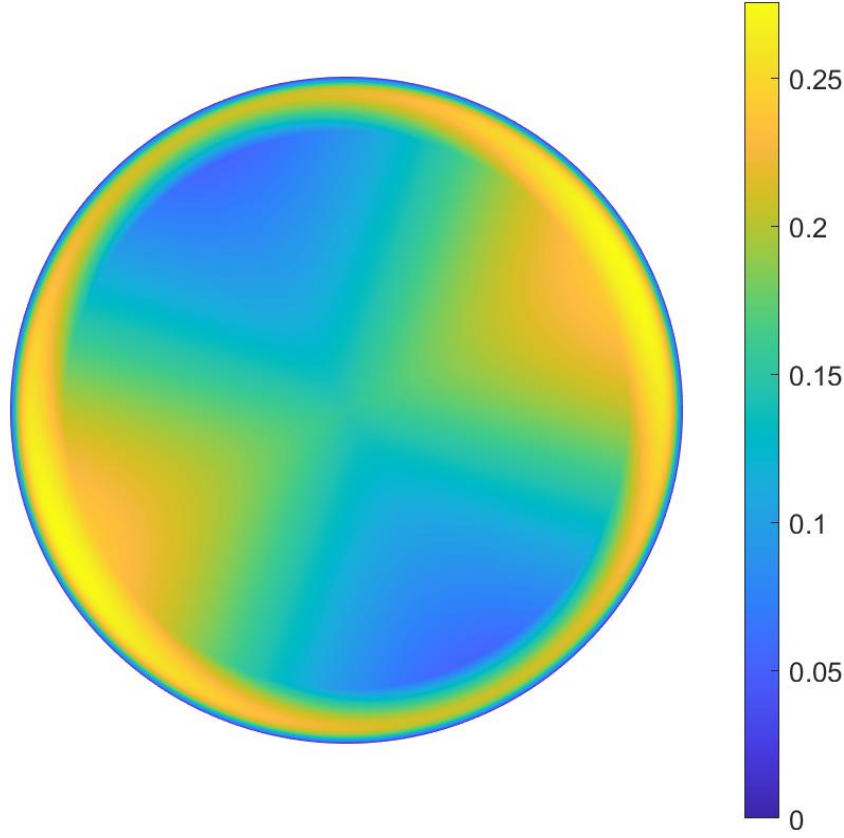

Figure 3.9 Mean errors distribution plot

## 3.4 Numerical Experiment 2

This section provides a more complex equation form and initial conditions. Additionally, numerical experiments were conducted in this part, and analyses and discussions were carried out based on the results of these numerical experiments. Consider the following equation form and initial conditions:

$$\frac{\partial \bm{n}}{\partial t} - (\frac{1}{r}\frac{\partial}{\partial r}\left(rD(r,\theta)\frac{\partial \bm{n}}{\partial r}\right) + \frac{1}{r^2}\frac{\partial}{\partial \theta}(D(r,\theta)\frac{\partial \bm{n}}{\partial \theta})) - (\frac{1}{r}\frac{\partial}{\partial r}(r\bm{n}\hat{\mu}(r,\theta)\bm{E}(r,\theta))$$

$$+\frac{1}{r}\frac{\partial}{\partial \theta}(\bm{n}\hat{\mu}(r,\theta)\bm{E}(r,\theta))) = GR(t,r,\theta), \quad 0<t<T, \quad (r,\theta) \in \Omega,$$

$$\bm{n}(t,\tilde{r},\theta) = e^{-t}\sin\theta^2 + \sin\theta^2, \quad 0<t<T, \quad (\tilde{r},\theta) \in \Gamma,$$

$$\bm{n}(0,r,\theta) = \sin\theta^2 r + r(1-r)^{1/2} + \sin\theta^2 r, \quad (r,\theta) \in \Omega,$$

$$\hat{\mu}(r,\theta)\bm{E}(r,\theta) = E(r,\theta) = r,$$

$$D(r,\theta) = r,$$

$$GR(t,r,\theta) = -e^{-t}\sin\theta^2 r - (2e^{-t}\sin\theta^2 + 2(1-r)^{1/2} - 2r(1-r)^{-1/2}$$

$$+\frac{1}{4}r^2(1-r)^{-3/2} + 2\sin\theta^2 + 2\cos(2\theta) + 2\cos(2\theta)e^{-t}) + 3e^{-t}r\sin\theta^2$$

$$+3r(1-r)^{1/2} - \frac{1}{2}r^2(1-r)^{-1/2} + 3r\sin\theta^2 + \sin 2\theta re^{-t} + \sin 2\theta r,$$

where $(r,\theta) \in \Omega = \{(r,\theta) \mid 0<r<1, 0 \leq \theta < 2\pi\}$ indicates that it is defined on a disk and has an exact solution $n(t,r,\theta) = e^{-t}\sin\theta^2 r + r(1-r)^{1/2} + \sin\theta^2 r$.

For $m > 0$, let $h = 1/(m+1)$ and $n = [2m\pi]$, which is the largest integer less than $2m\pi$, Take $r_i = \varphi(i*h)$ and $\tilde{\mu} = 2\pi/n$.

As shown in Table 3.2, when $t = 0.1$, $\Delta\tilde{t} = 0.01$, $a = 0.1$, $l = 5$, $p = 0.5$, $q = 5$, the results of grid division using three adaptive stretching functions and an adaptive difference method are compared. In this relatively complex example, the algorithm presented in this paper can improve accuracy across various dimensions. The improvement in algorithm accuracy is particularly significant in the radial, angular, and difference adaptive dimensions, with an overall accuracy increase of approximately 70%.

When $m = 99$ and $t = 0.1$, the exact solutions of the equation system are plotted against the approximate solutions when $m = 99$, $t = 0.1$, $\Delta\tilde{t} = 0.01$, $a = 0.1$, $l = 5$, $p = 0.5$, $q = 5$, as shown in Figures 3.10 and 3.11. Observing the three-dimensional plots of the exact and approximate solutions, it is evident that the shape of the approximate solution closely matches the exact solution, indicating good simulation results.

Figure 3.12 illustrates the average error variation with $h$ and $p$ when $t = 0.1$, $\Delta\tilde{t} = 0.01$, $a = 0.5$, $q = 1$, $l = 1$. The figure shows that the average error first decreases and then increases as $p$ increases, and it also indicates that the average error decreases as the grid becomes finer. This demonstrates that the algorithm is convergent. Furthermore, there exists an optimal value of $p$, approximately 0.5, which minimizes the average error.

Figure 3.13 displays the average error variation with $h$ and $q$ when $t = 0.1$, $\Delta\tilde{t} = 0.01$, $a = 0.5$, $p = 0.5$, $l = 1$. The figure shows that the average error decreases as $q$ increases, and it also indicates that the average error decreases as the grid becomes finer. This demonstrates that the algorithm is convergent. Additionally, there exists an optimal value of $q$, approximately 5, which minimizes the average error.

Table 3.2 Numerical experiment 2 Meanerr results comparison

| m | Meanerr(Uniform Adaptation) | Angle Non-adaptation | Radius Non-adaptation | Time Non-adaptation | Difference Non-adaptation | Non-adaptation |
|---|---|---|---|---|---|---|
| 99 | 1.97999984e-01 | 2.00826253e-01 | 3.06133431e-01 | 1.99343523e-01 | 2.37491614e-01 | 3.32342775e-01 |

Figure 3.14 shows the results of how the average error changes with $h$ and $l$ when $t = 0.1$, $\Delta\tilde{t} = 0.01$, $a = 0.5$, $p = 0.5$, $q = 5$. The figure indicates that the average error decreases consistently as $l$ increases, and it also shows that the average error decreases as the grid becomes finer. This demonstrates that the algorithm is convergent. Furthermore, there exists an optimal value of $l$, approximately 5, which minimizes the average error.

Figure 3.15 displays the results of how the average error changes with $h$ and $a$ when $t = 0.1$, $\Delta\tilde{t} = 0.01$, $l = 5$, $p = 0.5$, $q = 5$. The figure shows that the average error first increases and then decreases as $a$ increases, and it also indicates that there is an optimal value of $a$ around 0.1, which minimizes the average error when $t = 0.1$ and $\Delta \tilde{t} = 0.01$. Additionally, the figure demonstrates that the average error decreases as the grid becomes finer, indicating that the algorithm is convergent. Comparing the exploration of the optimal parameters in both experiments, in this relatively complex case, setting $a$ to a relatively small value, $p$ close to 0.5, and $q$ and $l$ to relatively large values can result in smaller errors. The experimental results also show that under different real-world conditions, the values of these four parameters significantly affect the results of numerical simulations. When these values are close to their optimal values, the grid in the region of intense oscillation becomes denser and contains more information, thereby improving the accuracy of numerical simulations compared to non-adaptive algorithms. The complexity of the algorithm presented in this paper does not change relative to non-adaptive algorithms; therefore, selecting better parameters based on actual conditions can greatly increase the accuracy of numerical simulations without changing the algorithm's efficiency.

At $t = 0.1$ and $\Delta \tilde{t} = 0.01$, the distribution of errors in the approximate solution under the conditions $h = 1/100$, $a = 0.1$, $l = 5$, $p = 0.1$, $q = 5$ is shown in Figure 3.16. The figure shows that the error distribution is relatively uniform. The angle stretching function and radial stretching function can significantly reduce the error at specific positions, thereby reducing the overall average error and improving the overall simulation effect.

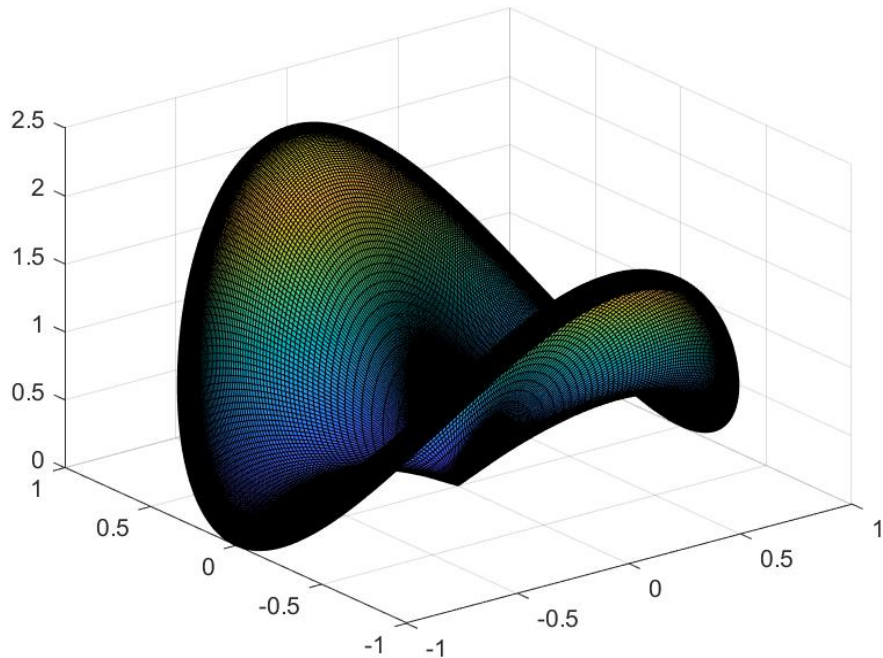

**Figure 3.10 Exact solution 3D plot**

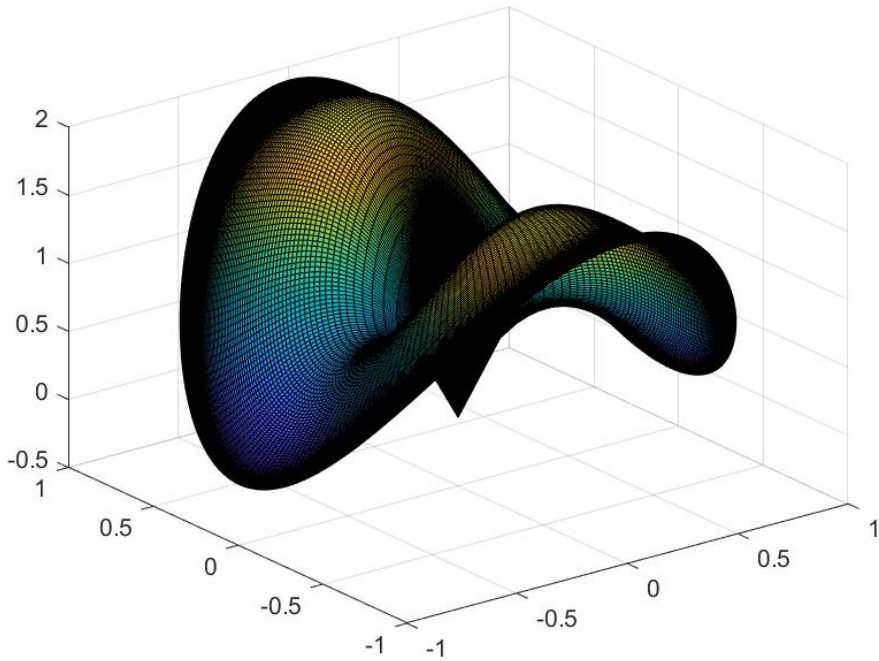

**Figure 3.11 Approximate solution 3D plot**

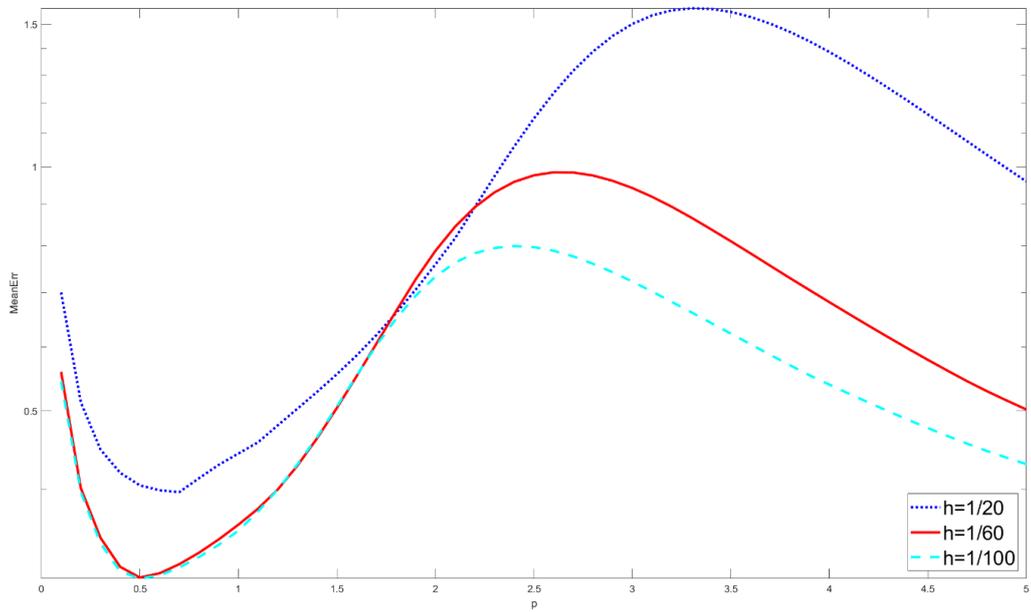

**Figure 3.12 Mean errors with respect to $p$ and $h$ plot**

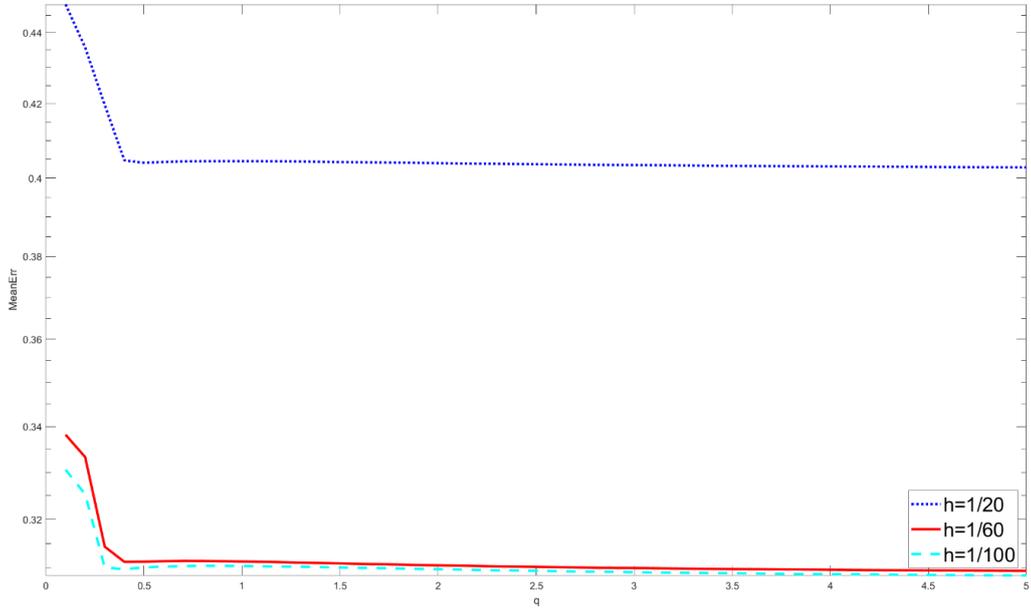

**Figure 3.13 Mean errors with respect to** $q$ **and** $h$ **plot**

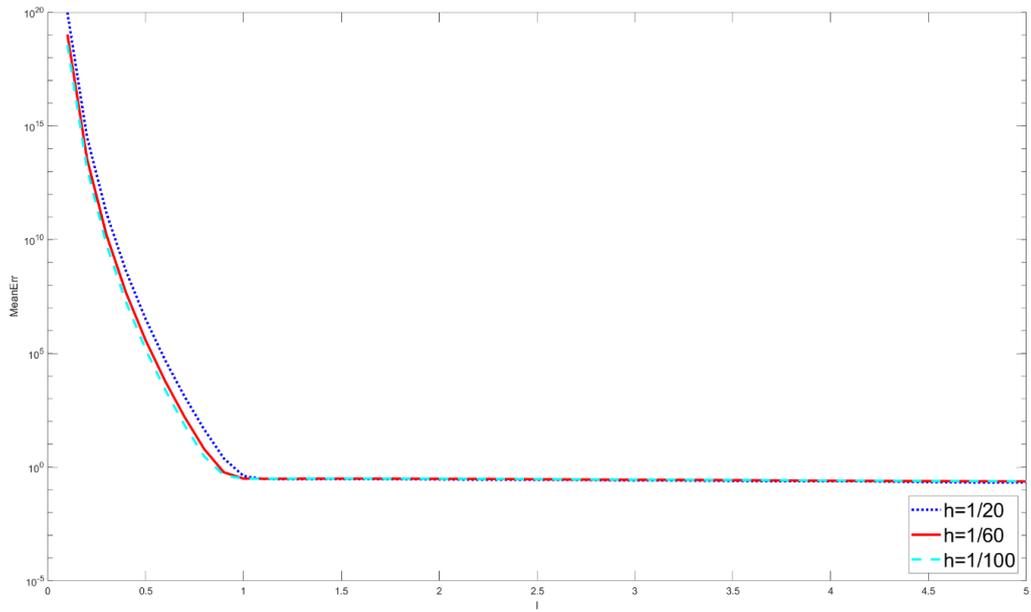

**Figure 3.14 Mean errors with respect to** $l$ **and** $h$ **plot**

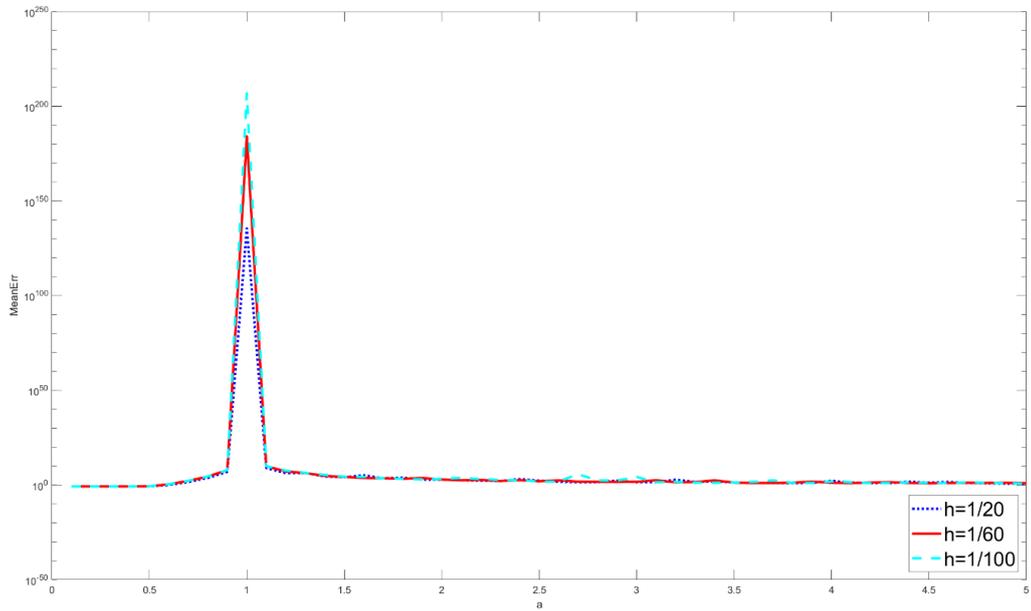

**Figure 3.15 Mean errors with respect to $a$ and $h$ plot**

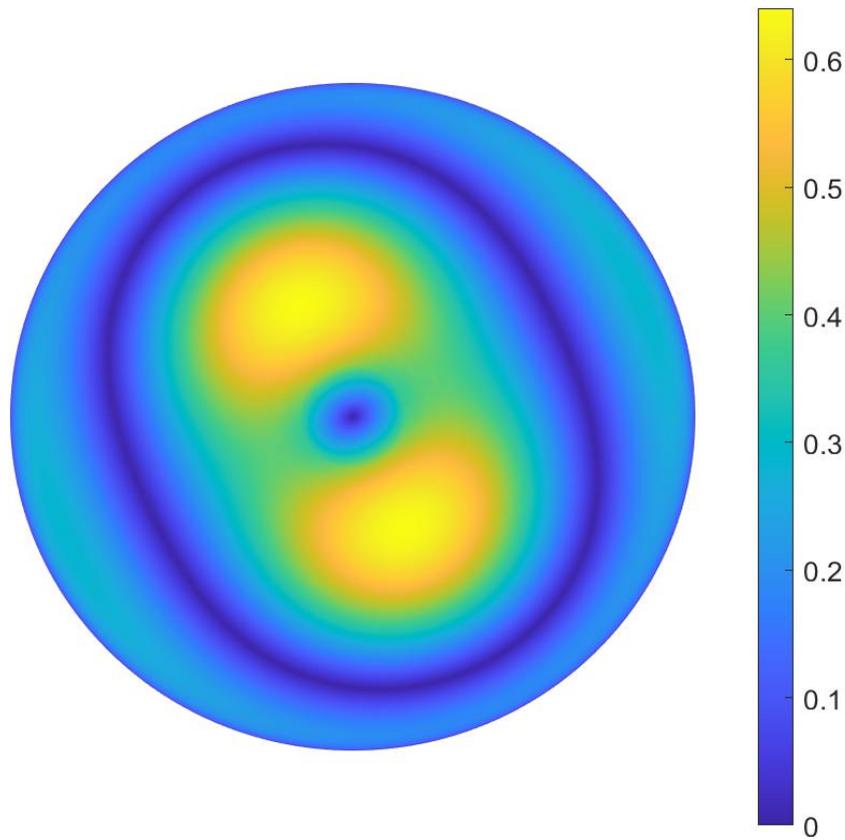

**Figure 3.16 Mean errors distribution plot**

# 4 Conclusion and Outlook

This chapter systematically summarizes the content of the previous three chapters and elaborates on the significance and future prospects of the multi-dimensional self-adaptive numerical simulation framework presented in this paper.

## 4.1 Conclusion

In Chapter One, the introduction elaborates on the significant roles of two important extended forms of the Boltzmann transport equation for semiconductor devices—the parabolic equation and the continuity equations for electrons and holes—across various fields. It highlights the limitations of traditional finite difference methods in handling regions with singular properties, nonlinear terms, and computational efficiency. The chapter systematically reviews the research background and challenges associated with numerical simulations using finite difference methods for semiconductors. To address existing issues in current numerical simulation studies, it proposes innovative research content on multi-dimensional self-adaptive numerical simulations, including coordinate transformation of equation systems into polar form, adaptive drift and diffusion coefficients, multi-dimensional self-adaptive mesh partitioning, and adjustable-parameter difference schemes. This underscores the theoretical and practical value of the present study.

In Chapter Two, the paper addresses the self-adaptive solution of a parabolic partial differential equation for semiconductors defined on a unit disk with boundary singularities. A radial-direction self-adaptive mesh partitioning method is proposed, combined with a Crank-Nicolson scheme incorporating tunable parameters. The Swartztrauber-Sweet method is employed to construct the finite difference formulation at the origin, effectively eliminating singularity-related issues arising from the pole in polar coordinates and boundary positions. Convergence analysis indicates that both maximum and average errors converge to the product of a constant and a function dependent on step size and parameters. Subsequent numerical experiments verify the convergence and superiority of the proposed method. In Example One, results show that the algorithm significantly improves numerical accuracy and enhances algorithmic stability (with computational efficiency unchanged, accuracy improves by more than sevenfold, and the algorithm converges). In the more complex Example Two, the algorithm again demonstrates improved accuracy and stability (accuracy increases approximately twofold under unchanged computational efficiency, and the algorithm converges). Furthermore, the study conducts sensitivity analyses of errors with respect to parameter variations and mesh refinement levels, along with visualizations of numerical solutions and error distributions. Optimal self-adaptive parameters (for radius and difference schemes) are determined, providing guidance for practical parameter tuning.

Chapter Three aims to solve the self-adaptive numerical simulation problem for the semiconductor continuity equations derived from the Boltzmann equation on a unit disk exhibiting singular behavior. The continuity equations for electrons and holes are transformed into polar coordinates. By combining multi-dimensional self-adaptive grids in the radial, angular, and temporal dimensions, and introducing an adaptive difference scheme with adjustable parameters, a multi-dimensional self-adaptive numerical simulation technique is constructed. The control volume integration method is applied to handle discretization at the origin. Treating drift and diffusion coefficients as functions varying with spatial position, this approach—combined with the polar coordinate representation—greatly enhances the interpretability of numerical simulations in real-world scenarios and their adaptability to different semiconductor device simulations. In the simpler case (Example One), the algorithm improves accuracy in the radial dimension and enhances algorithmic stability (accuracy increases by approximately 10% without sacrificing computational efficiency, and the algorithm converges). In the more complex Example Two, the algorithm achieves a significant improvement in numerical accuracy (approximately 70% increase under unchanged computational efficiency). Numerical experiments demonstrate that the algorithm substantially enhances solution accuracy and numerical stability for complex nonlinear problems, avoiding numerical oscillations caused by singularities. Compared to uniform mesh partitioning across all dimensions, the proposed method achieves much higher precision, ensures algorithmic convergence, exhibits better stability, and adaptively adjusts the mesh based on solution characteristics. Moreover, the optimal selection of adaptive parameters (such as angular parameters and weighting coefficients in difference schemes) has a significant impact on accuracy. This study further explores optimal parameter settings through comparative visualization of numerical results, verifies the convergence and superiority of the algorithm, and provides practical parameter-tuning strategies.

In summary, this research constructs a high-precision self-adaptive numerical simulation framework: coordinate transformation into polar form and variable drift and diffusion coefficients make the simulations more consistent with the physical characteristics of semiconductor devices; multi-dimensional self-adaptive grids and adjustable-parameter difference schemes are used to capture singular regions; the Swartztrauber-Sweet method and control volume integration eliminate the singularity at the origin induced by polar coordinate transformation; and parallelized code implementation ensures computational efficiency. The self-adaptive algorithm is applied to numerical experiments on both major extensions of the semiconductor Boltzmann transport equation, verifying its convergence, stability, and superiority, while also providing guidance for parameter tuning and simulation strategies. The outcomes offer methodological support and efficient tools for high self-adaptive numerical simulation and optimization design of semiconductor devices, addressing the issue of insufficient accuracy in traditional methods when solving semiconductor equations in singular regions. The use of parallelized MATLAB programming balances precision and computational efficiency, thereby providing theoretical support for high self-adaptivity and efficiency in semiconductor device simulations.。

## 4.2 Outlook

Although this study has achieved certain results in the field of adaptive numerical methods and semiconductor device simulation, there are still directions that require further exploration and improvement:

This research primarily focuses on numerical simulations of single physical fields, whereas actual semiconductor devices often involve coupled electro-thermal-mechanical effects. Furthermore, when designing nanoscale semiconductor devices, quantum effects within the semiconductor also play a significant role. To more accurately simulate the internal mechanisms of semiconductor devices, future work should couple the drift-diffusion equations with heat conduction equations and mechanical deformation equations, and incorporate density functional theory (DFT) or quantum correction terms, thereby constructing a more comprehensive multi-physics model.

The selection of adaptive parameters in this study mainly relies on empirical tuning, which is relatively inefficient. If the adaptive algorithm is integrated with machine learning to establish a self-adaptive decision-making model for parameters, neural networks could be used to predict local error distributions and guide real-time mesh refinement strategies, achieving dynamic optimal control. This would greatly enhance the adaptability and generalizability of the algorithm.

Although the numerical simulation program developed in this study has implemented multi-threaded parallelization, it still faces memory and computational bottlenecks under ultra-large-scale grids. In the future, GPU acceleration technologies (such as CUDA architecture) and distributed computing frameworks (such as MPI) should be combined to improve the scalability of the algorithm. Additionally, developing efficient preconditioners based on algebraic multigrid (AMG) methods could further optimize the solution efficiency of linear systems.

The algorithm presented in this paper is only applicable to static circular disk domains and cannot handle numerical simulations for devices with irregular geometries or dynamically deforming boundaries (such as MEMS devices). Future efforts should focus on developing immersed boundary methods (IBM) or moving mesh techniques, combined with adaptive strategies, to address dynamic interface problems.

Since this study is largely theoretical and has limited integration with practical applications, future work should strengthen collaboration with semiconductor manufacturing enterprises. Through experimental data (such as carrier mobility and interface barrier height), model parameters can be inverted to enhance prediction accuracy.

In summary, the multi-dimensional self-adaptive numerical simulation framework presented in this paper shows significant potential in the simulation of semiconductor devices. However, it still requires further refinement through interdisciplinary collaboration and technological innovation. With the advancement of computational technology and the refinement of physical models, high-precision and high-efficiency multi-physics coupled numerical simulation techniques will become the core tools for the next generation of semiconductor design. These advancements will facilitate independent innovation and breakthroughs in the field of integrated circuits.